\documentclass[lettersize,journal]{IEEEtran}
\usepackage{amsmath,amsfonts}
\usepackage{algorithmic}
\usepackage{algorithm}
\usepackage{array}
\usepackage[caption=false,font=normalsize,labelfont=sf,textfont=sf]{subfig}
\usepackage{textcomp}
\usepackage{stfloats}
\usepackage{url}
\usepackage{verbatim}
\usepackage{graphicx}
\usepackage{cite}
\usepackage{booktabs}
\usepackage[numbers]{natbib}
\usepackage{multirow}
\usepackage{multicol}
\usepackage{authblk}
\usepackage{diagbox}

\begin{document}

\title{Infrared Small Target Detection via tensor $L_{2,1}$ norm minimization and ASSTV regularization: A Novel Tensor Recovery Approach}
\author[a]{
Jiqian Zhao and An-Bao Xu$^{*}$\thanks{* Corresponding author}\thanks{This work was supported by the National Natural Science Foundation of
China under Grant 11801418.}\thanks{Jiqian Zhao is with the College of Computer Science and artificial intelligence, Wenzhou University, Zhejiang 325035, China (e-mail:zjq18698557926@gmail.com).}\thanks{An-Bao Xu is with the College of Mathematics and Physics, Wenzhou University, Zhejiang 325035, China (e-mail:xuanbao@wzu.edu.cn).}
}


\renewcommand*{\Affilfont}{\small\it}


\maketitle

\begin{abstract}
In recent years, there has been a noteworthy focus on infrared small target detection, given its vital importance in processing signals from infrared remote sensing. The considerable computational cost incurred by prior methods, relying excessively on nuclear norm for noise separation, necessitates the exploration of efficient alternatives. The aim of this research is to identify a swift and resilient tensor recovery method for the efficient extraction of infrared small targets from image sequences.
Theoretical validation indicates that smaller singular values predominantly contribute to constructing noise information. In the exclusion process, tensor QR decomposition is employed to reasonably reduce the size of the target tensor. Subsequently, we address a tensor $L_{2,1}$ Norm Minimization via T-QR (TLNMTQR) based method to effectively isolate the noise, markedly improving computational speed without compromising accuracy.
Concurrently, by integrating the asymmetric spatial-temporal total variation regularization method (ASSTV), our objective is to augment the flexibility and efficacy of our algorithm in handling time series data.
Ultimately, our method underwent rigorous testing with real-world data, affirmatively showcasing the superiority of our algorithm in terms of speed, precision, and robustness.

\end{abstract}


\section{INTRODUCTION}
\IEEEPARstart{I}{n} the contemporary era of the internet, our lives are intricately woven into the fabric of signals \cite{signals_matrix}. As we harness these signals, our constant aspiration is for precise and expeditious data collection. Nonetheless, in practical application, the dichotomy between speed and accuracy arises due to technological constraints. Sampling at lower rates presents a formidable challenge in acquiring truly comprehensive data. This challenge is particularly pronounced in the domain of infrared small target detection \cite{ISTD1, ISTD2, ISTD3}, where there is a pressing demand for a rapid and efficient methodology.

Compressive sensing \cite{compress_sensing, 8260873} offers a technical solution by sparsely sampling signals at lower rates and solving a optimization problems to ultimately reconstruct the complete signal. However, this method often applies to one-dimensional data. In contemporary society, data is typically stored in matrices, such as infrared images, audio, and video. Similarly, these data encounter partial element loss during collection, transmission, and storage. If compressive sensing methods are used, they may overlook the relationships between elements in two-dimensional rows and columns. 

This is where matrix recovery methods come into play. When discussing early effective algorithms in matrix recovery, it's inevitable to mention Robust Principal Component Analysis (RPCA). The RPCA model is defined as $\mathbf{Y} = \mathbf{L} + \mathbf{S}$. This model describes the process of decomposing a noisy original matrix $\mathbf{Y}$ into a low-rank matrix $\mathbf{L}$ and a sparse matrix $\mathbf{S}$. In the solving process, an optimization problem composed of the $L_1$-norm and nuclear norm is used. These methods fall under the category of Low-Rank and Sparse Matrix Decomposition Models (LRSD) \cite{CHANDRASEKARAN20091493, LRSD1, LRSD2}. For instance, \cite{lin2010augmented} introduced a swift and effective solver for robust principal component analysis. Building upon this, various methods have since emerged, leveraging the concept of low-rank matrix recovery, see \cite{chen2017denoising, xie2016hyperspectral, zhang2013hyperspectral}.

\begin{table*}[!htb]
	\caption{Explanation of notation in this paper}
	\centering
	\resizebox{\linewidth}{!}
	{
		\begin{tabular}{cc|cc}
			\toprule
			\toprule
			$\mathbf{a}$                                                      & vectors                                        & $\mathbf{I}_{n}$                                         & the identity matrix                               \\
			$\mathbf{A}$                                                      & matrices                                       & $\mathcal{I}$                                            & the identity tensor                               \\
			$\mathcal{A}$                                                     & tensor                                         & $\left \langle \mathcal{A},\mathcal{B}   \right \rangle$ & the inner product of $\mathcal{A},\mathcal{B}$    \\
			$\mathcal{A} \in \mathbb{C} ^{n_{1}\times n_{2}\times n_{3}   } $ & a third-order tensor                           & conj$ \left ( \mathcal{A}  \right ) $                    & the complex conjugate of $ \mathcal{A}$           \\
			$\mathcal{A}\left ( i,j,k \right ) $                              & the $(i,j,k)\mbox{-}$th entry of $\mathcal{A}$ & $\left \| \mathcal{A}  \right \| _{F}$                   & the tensor Frobenius norm of $\mathcal{A}$        \\
			$\mathcal{A}\left ( i,j, \right ) $                               & the tube fiber of $\mathcal{A}$                & $\left \| \mathcal{A}  \right \| _{1}$                   & the tensor $L_{1}\mbox{-}$norm of $\mathcal{A}$   \\
			$\mathcal{A}^{\left ( i \right ) } $                              & the $i\mbox{-}$th frontal slice $\mathcal{A}$  & $\left \| \mathcal{A}  \right \| _{2,1}$                 & the tensor $L_{2,1}\mbox{-}$norm of $\mathcal{A}$ \\
			$ \mathcal{A}^{\ast}$                                             & the conjugate transpose of $ \mathcal{A}$      & $\left \| \mathcal{A}  \right \| _{*}$                   & the tensor nuclear norm of $\mathcal{A}$          \\
			$\hat{\mathcal{A}}$                                               & the result of DFT on $\mathcal{A}$             & $\left \| \mathcal{A}  \right \| _{ASSTV}$               & the tensor ASSTV norm of $\mathcal{A}$            \\
			\bottomrule
			\bottomrule
		\end{tabular}
	}
	\label{table1}
\end{table*}

As research progresses, tensors are seen as a way to further optimize and upgrade LRSD, known as Tensor Low-Rank and Sparse Matrix Decomposition Models (TLRSD) \cite{TRPCA, TLRSD1}. Tensors, due to their increased dimensions, encompass more information and exhibit significantly enhanced computational efficiency. 

Analogous to matrix recovery methods, Tensor Robust Principal Component Analysis (TRPCA) \cite{8606166} was proposed and showed promising results in tensor recovery. Similarly, methods related to tensor decomposition have gradually emerged. In recent years, Total Variation (TV) regularization \cite{7892843, TV, TV2} has been widely applied to solve tensor recovery problems. However, traditional TV regularization, while considering spatial information, is computationally complex. Hence, spatial spectral total variation (SSTV) \cite{sstv} was introduced, significantly reducing the complexity of solving TV regularization. Yet, besides the spatial information within the tensor itself, the temporal information in tensors remained underutilized, leading to the introduction of ASSTV \cite{ASSTV} to incorporate both temporal and spatial information.

Thus, tensor recovery has evolved to primarily solve the ASSTV regularization and nuclear norm optimization problems. However, this process involves tensor singular value decomposition (TSVD) \cite{TSVD} for solving tensor nuclear norm minimization problem, which is highly time-consuming. To address this, a new method ASSTV-TLNMTQR was proposed, which combining ASSTV regularization and TLNMTQR. Our method contributes in the following ways:

\begin{enumerate}
    \item Introducing the innovative ASSTV-TLNMTQR method which has significantly elevated the pace of tensor decomposition. Through rigorous experimentation in infrared small target detection, this approach has showcased its rapid and efficient tensor recovery capabilities.

    \item A pioneering strategy, rooted in tensor decomposition, has been presented for minimizing the $L_{2,1}$ norm. This method has proven to be highly effective in solving tensor decomposition problems, surpassing the speed of previous approaches reliant on nuclear norm methodologies.

    \item Our approach integrates ASSTV regularization and TLNMTQR, resulting in a model endowed. This integration empowers our model to produce outstanding solutions across varied contexts, facilitated by flexible parameter adjustments.
\end{enumerate}

\section{NOTATIONS AND PRELIMINARIES}
\subsection{Fast fourier transform}
In this paper, we have compiled a comprehensive list of symbols and definitions essential for reference, as presented in TABLE \ref{table1}. In order to effectively articulate our model, it is imperative to introduce fundamental definitions and theorems. Primarily, we emphasize a critically significant operation associated with tensor product and matrix product— the Discrete Fourier Transform (DFT). Here, we denote $\hat{\mathbf{v}}$ as an integral part of this discussion.
\begin{center}
	$\hat{\mathbf{v}} = \mathbf{F_{n}\mathbf{v}}$
\end{center}
where $\mathbf{F_{n}}$ is DFT matrix denoted as
\begin{center}
	$\mathbf{F} _{n} = \begin{bmatrix}
			1      & 1             & 1                               & \cdots & 1                                                  \\
			1      & \omega        & \omega ^{2}                     & \cdots & \omega ^{n-1}                                      \\
			\vdots & \vdots        & \vdots                          & \ddots & \vdots                                             \\
			1      & \omega ^{n-1} & \omega ^{2\left ( n-1 \right )} & \cdots & \omega ^{\left ( n-1 \right )\left ( n-1 \right )}
		\end{bmatrix} \in \mathbb{C} ^{n\times n} $,
\end{center}
where $\omega =e^{-\frac{2\pi i}{n} }$ and $i$ is the imaginary unit. So we can learn that $\frac{\mathbf{F} _{n}}{\sqrt{n}}$ is orthogonal, i.e.,
\begin{equation}
	\label{e2.1.1}
	\mathbf{F}_{n}^{\ast}\mathbf{F}_{n}=\mathbf{F}_{n}\mathbf{F}_{n}^{\ast}=n\mathbf{I}_{n}.
\end{equation}
By providing a proof, we can easily obtain the following conclusions:
\begin{equation}\label{2.2}
	\left ( \mathbf{F} _{{n}_{3}}\otimes \mathbf{I}_{{n}_{1}}  \right ) \cdot \text{bcirc}\left ( \mathcal{A}  \right )\cdot \left ( \mathbf{F} _{{n}_{3}}^{-1}\otimes \mathbf{I}_{{n}_{2}} \right )  =\bar{\mathbf{A}  } ,
\end{equation}
where $\otimes$ denotes the Kronecker product and $\frac{ \mathbf{F} _{{n}_{3}}^{-1}\otimes \mathbf{I}_{{n}_{2}}}{\sqrt{n_{3}}}$ is orthogonal. Moreover, $\bar{\mathbf{A} }$ and $\text{bcirc}\left (\mathcal{A}\right )$ are defined as follows.
\begin{center}
	$\bar{\mathbf{A} } = \text{bdiag}\left(\hat{\mathcal{A}}\right) = \begin{bmatrix}
			\hat{\mathcal{A}}^{\left ( 1 \right ) } &                                         &        &                                             \\
			                                        & \hat{\mathcal{A}}^{\left ( 2 \right ) } &        &                                             \\
			                                        &                                         & \ddots &                                             \\
			                                        &                                         &        & \hat{\mathcal{A}}^{\left ( n_{3} \right ) }
		\end{bmatrix}$,
\end{center}
\begin{center}
	$\text{bcirc}\left (\mathcal{A}\right ) = \begin{bmatrix}
			\mathcal{A}^{\left ( 1 \right )}     & \mathcal{A}^{\left ( n_{3} \right )}   & \cdots & \mathcal{A}^{\left ( 2 \right )} \\
			\mathcal{A}^{\left ( 2 \right )}     & \mathcal{A}^{\left ( 1 \right )}       & \cdots & \mathcal{A}^{\left ( 3 \right )} \\
			\vdots                               & \ddots                                 & \ddots & \vdots                           \\
			\mathcal{A}^{\left ( n_{3} \right )} & \mathcal{A}^{\left ( n_{3}-1 \right )} & \cdots & \mathcal{A}^{\left ( 1 \right )}
		\end{bmatrix}$	.
\end{center}
\subsection{Tensor product}
Here, we provide the definition of the \textbf{T-product}:
Let $\mathcal{A}\in \mathbb{R}^{n_{1}\times n_{2}\times n_{3}}$ and $\mathcal{B}\in \mathbb{R}^{n_{2}\times l\times n_{3}}$. Then the t-product $\mathcal{A}\ast\mathcal{B}$ is the following tensor of size $n_{1}\times l\times n_{3}$:
\begin{equation}
	\mathcal{A}\ast\mathcal{B}=\text{fold}\left(\text{bcirc}\left(\mathcal{A}\right)\cdot\text{unfold}\left(\mathcal{B}\right)\right),
\end{equation}
where $\text{unfold}\left( \cdot \right)$ and $\text{fold}\left( \cdot \right)$ are the unfold operator map and fold operator map of $\mathcal{A}$ respectively, i.e.,
\begin{center}
    $\text{unfold}\left(\mathcal{A}\right)=\begin{bmatrix}
		\mathcal{A}^{\left ( 1 \right ) }  \\
		\mathcal{A}^{\left ( 2 \right ) } \\
		\vdots  \\
		\mathcal{A}^{\left ( n_{3} \right ) } 
	\end{bmatrix}\in \mathbb{R} ^{n_{1}n_{3}\times n_{2}}$, \\$\text{fold}\left(\text{unfold}\left(\mathcal{A}\right)\right)=\mathcal{A}$.
\end{center}

\subsection{Tensor QR decomposition}
The most important component in $L_{2,1}$-norm \cite{ZHENG2021108240} is the tensor QR decomposition, which is an approximation of SVD.
Let $\mathcal{A}\in \mathbb{R}^{n_{1}\times n_{2}\times n_{3}}$. Then it can be factored as
\begin{equation}
	\mathcal{A}=\mathcal{Q}\ast\mathcal{R},
\end{equation}
where $\mathcal{Q}\in \mathbb{R}^{n_{1}\times n_{1}\times n_{3}}$ is orthogonal, and $\mathcal{R}\in \mathbb{R}^{n_{1}\times n_{2}\times n_{3}}$ is analogous to the upper triangular matrix.

\subsection{Asymmetric spatialtemporal total variation}

The TV regularizaiton \cite{7892843, TV, TV2} is widely applied in infrared small target detection, as it allows us to better utilize the information in the images.By incorporating the TV regularizaiton, we can maintain the smoothness of the image during the processes of image recovery and denoising, thereby eliminating potential artifacts. Typically, the TV regularizaiton is applied to matrices as shown in Equation (\ref{TV}).

\begin{equation}\label{TV}
	{\left\| {\cal A} \right\|_{TV}} = {\left\| {{D_h}{\cal A}} \right\|_1} + {\left\| {{D_v}{\cal A}} \right\|_1}
\end{equation}

However, this approach neglects the temporal information. Hence, we use a new method called asymmetric spatial-temporal total variation regularizaiton method to model both spatial and temporal continuity.
There are two main reasons for choosing ASSTV regularizaiton method:
\begin{itemize}
	\item It facilitates the preservation of image details and better noise elimination.
	\item By introducing a new parameter $\delta$, it allows us to assign different weights to the temporal and spatial components, enabling more flexible tensor recovery.
\end{itemize}
Now, we provide the definition of ASSTV, as shown in Equation (\ref{ASSTV}):
\begin{equation}\label{ASSTV}
	{\left\| \mathcal{A} \right\|_{ASSTV}} = {\left\| {{\rm{ }}{D_h}\mathcal{A}} \right\|_1} + {\left\| {{\rm{ }}{D_v}\mathcal{A}} \right\|_1} + \delta {\left\| {{\rm{ }}{D_z}\mathcal{A}} \right\|_1}
\end{equation}
Where $D_h$, $D_v$, and $D_z$ are the horizontal, vertical, and temporal difference operators, respectively. $\delta$ denotes a positive constant, which is used to control the contribution in the temporal dimension. Here are the definitions of the three operators:
\begin{equation}
	\mathop D\nolimits_h \mathcal{A}(i,j,k) = \mathcal{A}(i + 1,j,k) - \mathcal{A}(i,j,k)
\end{equation}
\begin{equation}
	\mathop D\nolimits_v \mathcal{A}(i,j,k) = \mathcal{A}(i,j + 1,k) - \mathcal{A}(i,j,k)
\end{equation}
\begin{equation}
	\mathop D\nolimits_z \mathcal{A}(i,j,k) = \mathcal{A}(i,j,k + 1) - \mathcal{A}(i,j,k)
\end{equation}

\subsection{Tensor $L_{2,1}$ norm minimization via tensor QR decomposition}

\begin{figure}[h]
    \centering
        \caption{Singular value of one image}
    \includegraphics[width=0.9\linewidth]{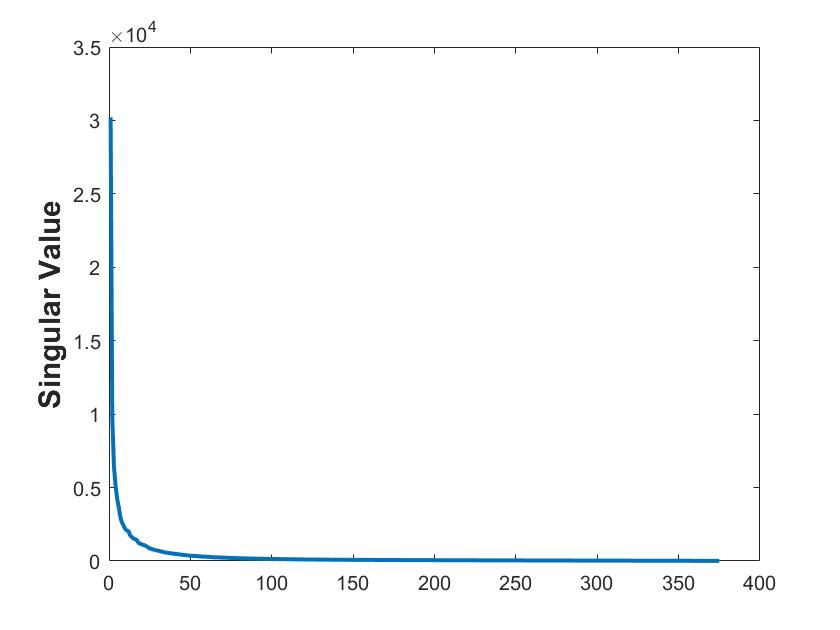}

    \label{sv}
\end{figure} 

In the realm of matrix recovery, a commonly employed approach involves solving the nuclear norm minimization problem through Singular Value Decomposition (SVD) to disentangle noise and restore images. When applying the SVD to a single channel of an RGB image, the singular values can be observed, as depicted in Fig. \ref{sv}. However, the SVD is computationally demanding and often requires a substantial amount of algorithmic execution time to converge. Consequently, we adopt a novel tensor recovery method based on tensor tri-factorization to address the \(L_{2,1}\) norm minimization problem. This approach aims to achieve image recovery with reduced computational complexity.
\begin{equation}
	\label{q1}
	\min \tau {\left\| \mathbf{Z} \right\|_{2,1}} + \frac{1}{2}\left\| {\mathbf{Z} - \mathbf{Y}} \right\|_F^2
\end{equation}

First, we will introduce the matrix version. For the problem (\ref{q1}), there exists an optimal solution (\ref{s1}).
\begin{equation}
	\label{s1}
	\mathbf{Z}\left( {:,j} \right) = \max \left\{ {\frac{{{{\left\| {\mathbf{Y}\left( {:,j} \right)} \right\|}_F} - \tau }}{{{{\left\| {\mathbf{Y}\left( {:,j} \right)} \right\|}_F}}},0} \right\}\mathbf{L}\left( {:,j} \right)
\end{equation}

However, computing the entire matrix \(Z\) requires significant computational resources and time. Moreover, in most cases, smaller singular values in the matrix are primarily used to reconstruct noise information within the matrix. Therefore, it is unnecessary to compute them during the process of matrix recovery. Hence, it aims to find a method that only computes the first few larger singular values and their corresponding singular vectors, accelerating the overall computation process.

This problem is essentially finding an approximate rank-\(r\) approximation of the matrix \(\mathbf{Z}\), meaning solving the problem (\ref{q1.5}).
\begin{equation}
	\label{q1.5}
	\mathop {\min }\limits_{\mathbf{L},\mathbf{D},\mathbf{R}} \left\| {\mathbf{Z} - \mathbf{L}\mathbf{D}\mathbf{R}} \right\|_F^2 \le \varepsilon
\end{equation}

Where \(\mathbf{L} \in \mathbb{R}^{m \times r}\), \(\mathbf{D} \in \mathbb{R}^{r \times r}\), \(\mathbf{R} \in \mathbb{R}^{r \times n}\), and \(r\) represents the rank of matrix \(\mathbf{Z}\) (\(r \in (0, n]\)). \(\mathbf{L}\) and \(\mathbf{R}\) are column orthogonal and row orthogonal matrices respectively, while \(\mathbf{D}\) is a regular square matrix. \(\varepsilon\) is a positive error threshold. For problem (\ref{q1.5}), we can use the iterative method to solve it.

We initialize three matrices as follows: \(\mathbf{L}_1 = \text{eye}(m, r)\), \(\mathbf{D}_1 = \text{eye}(r, r)\), and \(\mathbf{R}_1 = \text{eye}(r, n)\). Then, we iterate to update \(\mathbf{L}\), \(\mathbf{D}\), and \(\mathbf{R}\) using the QR decomposition. First, fixing \(\mathbf{R}\), we have:
\[ \left[ {\mathbf{L}_{j + 1}, \sim } \right] = \rm{qr}(\mathbf{Z}\mathbf{R}_{j}^T).\]
Next, fixing \(\mathbf{L}\), we obtain:
\[ \left[ {\mathbf{R}_{j + 1}, \mathbf{L}_{j + 1}^T\mathbf{Z}\mathbf{R}_{j + 1}^T} \right] = \rm{qr}(\mathbf{Z}\mathbf{L}_{j + 1}).\]
Finally:
\[ {\mathbf{D}_{j + 1}} = \mathbf{L}_{j + 1}^\mathbf{T}\mathbf{Z}\mathbf{R}_{j + 1}^T.\]
After multiple iterations, the tri-factorization of \(\mathbf{Z}\) is completed.

Next, we utilize the obtained tri-factorization matrices to optimize our problem. First, we have \(\mathbf{L}\mathbf{D}\mathbf{R} \approx Z\), and then we establish the following \(L_{2,1}\) minimization model:

\begin{equation}
	\label{q2}
	\min \tau {\left\| {\mathbf{L}\mathbf{D}\mathbf{R}} \right\|_{2,1}}  + \frac{1}{2}\left\| {\mathbf{L}\mathbf{D}\mathbf{R} - \mathbf{Y}} \right\|_F^2
\end{equation}

Since we know that both \(\mathbf{L}\) and \(\mathbf{R}\) are orthogonal, Equation (\ref{q2}) simplifies to Equation (\ref{q3}).
\begin{equation}
	\label{q3}
	\begin{gathered}
		\min \tau {\left\| \mathbf{D} \right\|_{2,1}} + \frac{1}{2}\left\| {\mathbf{D} - {\mathbf{L}^T}\mathbf{Y}{\mathbf{R}^T}} \right\|_F^2
	\end{gathered}
\end{equation}

Let \(\mathbf{C} = \mathbf{L}^T\mathbf{Y}\mathbf{R}^T\), then we have Equation (\ref{q4}).
\begin{equation}
	\label{q4}
	\begin{gathered}
		\min \tau {\left\| \mathbf{D} \right\|_{2,1}} + \frac{1}{2}\left\| {\mathbf{D} - \mathbf{C}} \right\|_F^2 \\
	\end{gathered}
\end{equation}

We observe that this problem, resembling (\ref{q1}), can be solved using the contraction operator from (\ref{s1}), as shown in Equation (\ref{s2}).
\begin{equation}
	\label{s2}
	\mathbf{D}\left( {:,j} \right) = \max \left\{ {\frac{{{{\left\| {\mathbf{C}\left( {:,j} \right)} \right\|}_F} - \tau }}{{{{\left\| {\mathbf{C}\left( {:,j} \right)} \right\|}_F}}},0} \right\}\mathbf{C}\left( {:,j} \right)
\end{equation}

Similarly, for the tensor \({\cal Z}\), we need to solve the following optimization problem:
\begin{equation}
	\label{tq1}
	\mathop {\min } {\left\| {\cal Z} \right\|_{2,1}} + \frac{\mu }{2}\left\| {{\cal Z} - \left( {{\cal B} - \frac{y}{\mu }} \right)} \right\|_F^2
\end{equation}

We need to perform tri-factorization on the tensor \({\cal Z}\) and find three tensors such that \({{\cal L}{\cal D}{\cal R} \approx \cal Z} \), where \({\cal L} \in \mathbb{R}^{n_1 \times r \times n_3}\) and \({\cal R} \in \mathbb{R}^{r \times n_2 \times n_3}\) are orthogonal, and \({\cal D} \in \mathbb{R}^{r \times r \times n_3}\) is a tensor.

We obtain \({\cal L}\), \({\cal D}\), and \({\cal R}\) by solving the following optimization problem:
\begin{equation}
	\left\| {{\cal Z} - {\cal L}{\cal D}{\cal R}} \right\|_F^2 \le \varepsilon
\end{equation}

Then, we can utilize iterative method to obtain \({\cal L}\), \({\cal D}\), and \({\cal R}\). The final solution to the problem (\ref{tq1}) is transformed into solving the following problem:
\begin{equation}
	\label{tq2}
	\mathop {\min } {\left\| {{\cal L}{\cal D}{\cal R}} \right\|_{2,1}} + \frac{\mu }{2}\left\| {{\cal L}{\cal D}{\cal R} - \left( {{\cal B} - \frac{y}{\mu }} \right)} \right\|_F^2
\end{equation}

Simplifying (\ref{tq2}) leads to Equation (\ref{tq3}).
\begin{equation}
	\label{tq3}
	\mathop {\min } {\left\| {\cal D} \right\|_{2,1}} + \frac{\mu }{2}\left\| {{\cal D} - {{\cal L}^*}\left( {{\cal B} - \frac{y}{\mu }} \right){{\cal R}^*}} \right\|_F^2
\end{equation}

Let \({{{\cal L}^*}\left( {{\cal B} - \frac{y}{\mu }} \right){{\cal R}^*}} = {{\cal D}_T}\). The optimal solution for this modified problem is:
\begin{equation}
	{\cal D} = {\rm{ifft}}\left( {\max \left\{ {\frac{{{{\left\| {{{\cal D}_T}} \right\|}_F} - \tau}}{{{{\left\| {{{\cal D}_T}} \right\|}_F}}},0} \right\}{{\cal D}_T},[],3} \right)
\end{equation}

The complete algorithmic process is illustrated in Algorithm \ref{algorithm}.
\begin{algorithm}[H]
	\caption{: Tensor $L_{2,1}$-Norm Minimization Method Based on Tensor QR Decomposition.}
	\label{algorithm}
	\begin{algorithmic}[1]
		\REQUIRE
		A real incomplete tensor $\mathcal{X}\in \mathbb{R}^{n_{1}\times n_{2}\times n_{3}}$;
		\STATE Initialize: $r>0$,\quad $t>0$,\quad $\tau>0$,\quad $\varepsilon=1e-6$  $\mathcal{L}_{0}=\mathcal{I}\in \mathbb{R}^{n_{1}\times r\times n_{3}}$,\quad $\mathcal{D}_{0}=\mathcal{I}\in \mathbb{R}^{r\times r\times n_{3}}$,\quad $\mathcal{R}_{0}=\mathcal{I}\in \mathbb{R}^{r\times n_{2}\times n_{3}}$,\quad $\mathcal{M}=\mathcal{X}$.
            \WHILE{Not converged}
		\STATE $\left [ \mathcal{L},\sim    \right ] =\text{T-QR}\left(\mathcal{M} \ast \mathcal{R}^{\ast } \right)$;
		\STATE $\left [ \mathcal{R},\mathcal{D}_{T}^{\ast}    \right ] =\text{T-QR}\left(\left(\mathcal{M}\right)^{\ast}\ast \mathcal{L} \right)$;
		\STATE ${{\cal R}} = {\cal R}^*$
		\FOR{$t=1,...,n_{3}$}
		\FOR{$j=1,...,r$}
		\STATE ${\mathcal{D}^{j}}^{\left(t\right)}=\text{ifft}\left(\frac{\max \left \{ \left \| \hat{ \mathcal{D}_{T}^{j}}^{\left(t\right)} \right \|_{F} -\tau,0 \right \}  }{\left \| \hat{ \mathcal{D}_{T}^{j}}^{\left(t\right)} \right \|_{F}}  \hat{ \mathcal{D}_{T}^{j}}^{\left(t\right)} ,[],3\right) $;
		\ENDFOR
		\ENDFOR
            \STATE $\cal M = {\cal L}*{\cal D}*{\cal R} $
		\STATE Check the convergence conditions \\ 
            $\left\| {{\cal X} - {\cal L}{\cal D}{\cal R}} \right\|_F^2 \le \varepsilon$
		\ENDWHILE
		 \RETURN $\mathcal{M}$.
	\end{algorithmic}
\end{algorithm}

\section{OUR METHOD FOR TENSOR RECOVERY}
\begin{figure}[htbp]
\centering
\centering
\caption{The construction of the tensor}
\includegraphics[width=3.5cm]{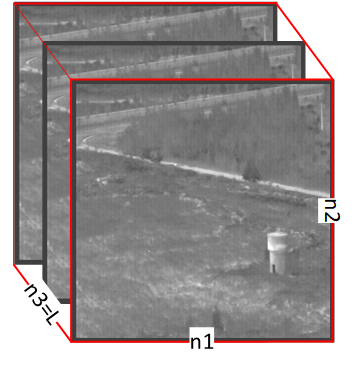}
\label{t1}
\end{figure}

There are different methods for constructing tensor $\cal K$ under different backgrounds. In infrared small target detection, we divide a continuous sequence of images into blocks by using a sliding window \cite{ASSTV, 6595533}, moving from the top left to the bottom right. Eventually, we obtain a tensor $\mathcal{K}$, see Fig. \ref{t1}. The height and width of the sliding window are $n_1$ and $n_2$, respectively, with a thickness of $n_3 = L$.

After constructing the tensor $\cal K$, if we performing the SVD on $\mathcal{K}$, we can clearly observe that $\mathcal{K}$ is a low-rank tensor.

Then, we define the following linear model:
\begin{equation}
	\mathcal{K} = \mathcal{B} + \mathcal{T} + \mathcal{N}
\end{equation}

Where ${\cal B}$, ${\cal T}$, and ${\cal N}$ represent the background image, target image, and noise image,respectively, ${\cal K},{\cal B},{\cal T},{\cal N} \in {^{{n_1} \times {n_2} \times {n_3}}}$.

	\subsection{The proposed model}

	By combining the $L_{2,1}$ norm and ASSTV regularization, we present the definition of our model as follows:
	\begin{equation}\label{nm}
		\begin{split}
			\mathcal{B},\mathcal{T},\mathcal{N} = \mathop {\arg \min }\limits_{{\cal B},{\cal T},{\cal N}} {\left\| \mathcal{B} \right\|_{2,1}} &+ {\lambda _{tv}}{\left\| \mathcal{B} \right\|_{ASTTV}} \\
			+ {\lambda _s}{\left\| \mathcal{T} \right\|_1} &+ {\lambda _3}\left\| \mathcal{N} \right\|_F^2 \\
		\end{split}
	\end{equation}
	\begin{center}
		$s.t. \quad \mathcal{K} = \mathcal{B} + \mathcal{T} + \mathcal{N},$
	\end{center}

	By utilizing the ASSTV formulation (\ref{ASSTV}), we can rewrite (\ref{nm}) as follows:
	\begin{equation}
		\begin{split}
			{\cal B},{\cal T},{\cal N} &= \mathop {\arg \min }\limits_{{\cal B},{\cal T},{\cal N}} {\left\| {\cal B} \right\|_{2,1}} + {\lambda _s}{\left\| {\cal T} \right\|_1} + {\lambda _3}\left\| {\cal N} \right\|_F^2\\
			&+ {\lambda _{tv}}({\left\| {{K_h}{\cal B}} \right\|_1} + {\left\| {{K_v}{\cal B}} \right\|_1} + \delta {\left\| {{K_z}{\cal B}} \right\|_1})
		\end{split}
	\end{equation}
	\begin{center}
		$s.t.\quad \mathcal{K} = \mathcal{B} + \mathcal{T} + \mathcal{N},$
	\end{center}

	Our model employs the $L_{2,1}$ norm, which significantly enhances algorithm efficiency. Additionally, we utilize tensor low-rank approximation to better assess the background, and the use of ASSTV allows for more flexible recovery of tensor.

	\subsection{Optimization procedure}
	In this section, we utilize the Alternating Direction Method of Multipliers (ADMM) method to solve (\ref{nm}). First, we modify (\ref{nm}) to the following form:
	\begin{equation}\label{cnm}
		\begin{split}
			{\cal B},{\cal T},{\cal N} &= \mathop {\arg \min }\limits_{B,T,N} {\left\| {\cal Z} \right\|_{2,1}} + {\lambda _s}{\left\| {\cal T} \right\|_1} + {\lambda _3}\left\| {\cal N} \right\|_F^2\\
			&+ {\lambda _{tv}}({\left\| {{{\cal V}_1}} \right\|_1} + {\left\| {{{\cal V}_2}} \right\|_1} + \delta{\left\| {{{\cal V}_3}} \right\|_1})
		\end{split}
	\end{equation}
	\begin{center}
		$s.t.\quad{\cal K} = {\cal B} + {\cal T} + {\cal N}, {\cal Z} = {\cal B},{{\cal V}_1} = {{\cal K}_h}{\cal B}, {{\cal V}_2} = {{\cal K}_v}{\cal B}, {{\cal V}_3} = {{\cal K}_z}{\cal B}$
	\end{center}
	Next, we present the augmented Lagrangian formulation of (\ref{cnm}):
	\begin{equation}
		\label{L}
		\begin{split}
			&{{L}_A}({\cal B},{\cal T},{\cal N},{\cal Z},{\cal V}) \\
			&= {\left\| {\cal Z} \right\|_{2,1}} + {\lambda _s}{\left\| {\cal T} \right\|_1} + {\lambda _3}\left\| {\cal N} \right\|_F^2 \\
			&+ {\lambda _{tv}}({\left\| {{{\cal V}_1}} \right\|_1} + {\left\| {{{\cal V}_2}} \right\|_1} + \delta{\left\| {{{\cal V}_3}} \right\|_1})\\
			&+ \left\langle {{y_1},{\cal K} - {\cal B} - {\cal T} - {\cal N}} \right\rangle  + \left\langle {{y_2},{\cal Z} - {\cal B}} \right\rangle  + \left\langle {{y_3},{{\cal V}_1} - {{\cal K}_h}{\cal B}} \right\rangle  \\
			&+ \left\langle {{y_4},{{\cal V}_2} - {{\cal K}_v}{\cal B}} \right\rangle  + \left\langle {{y_5},{{\cal V}_3} - {{\cal K}_z}{\cal B}} \right\rangle \\
			&+ \frac{\mu }{2}\left( {\left\| {{\cal K} - {\cal B} - {\cal T} - {\cal N}} \right\|_F^2} \right) + \left\| {{\cal Z} - {\cal B}} \right\|_F^2 \\
			&+ \left\| {{{\cal V}_1} - {{\cal K}_h}{\cal B}} \right\|_F^2 + \left\| {{{\cal V}_2} - {{\cal K}_v}{\cal B}} \right\|_F^2 + \left\| {{{\cal V}_2} - {{\cal K}_v}{\cal B}} \right\|_F^2
		\end{split}
	\end{equation}
	where $\mu$ represent a positive penalty scalar and $y_{1}$, $y_{2}$, $y_{3}$, $y_{4}$, $y_{5}$ represent the Lagrangian multiplier. We will use ADMM algorithm to solve (\ref{L}), which includes $\cal Z$, $\cal B$, $\cal T$, $V_{1}$, $V_{2}$, $V_{3}$ and $\cal N$. Then we will alternately update the variable as:
	\begin{enumerate}
		\item We update the value of ${\cal Z}$ using the following equation:
		      \begin{equation}\label{z}
			      {{\cal Z}^{k + 1}} = \mathop {\arg \min }\limits_Z {\left\| {\cal Z} \right\|_{2,1}} + \frac{{{\mu ^k}}}{2}\left\| {{\cal Z} - {\cal B} + \frac{{y_2^k}}{{{\mu ^k}}}} \right\|_F^2
		      \end{equation}
		      For more details about the update of ${\cal Z}$, please refer to Algorithm \ref{algorithm}. In the experiment, one iteration produced promising results, so we only performed one.
		\item We update the value of ${\cal B}$ using the following equation:
		      \begin{equation}\label{b}
			      \begin{array}{l}
				      {{\cal B}^{k + 1}} = \frac{{{\mu ^k}}}{2}\left( {\left\| {{\cal K} - {\cal B} - {{\cal T}^k} - {{\cal N}^k} + \frac{{y_1^k}}{{{\mu ^k}}}} \right\|_F^2} \right.               \\
				      + \left\| {{{\cal Z}^{k + 1}} - {\cal B} + \frac{{y_2^k}}{{{\mu ^k}}}} \right\|_F^2 + \left\| {{\cal V}_1^k - {{\cal K}_h}{\cal B} + \frac{{y_3^k}}{{{\mu ^k}}}} \right\|_F^2 \\
				      \left. { + \left\| {{\cal V}_2^k - {{\cal K}_v}{\cal B} + \frac{{y_4^k}}{{{\mu ^k}}}} \right\|_F^2 + \left\| {{\cal V}_3^k - {{\cal K}_z}{\cal B} + \frac{{y_5^k}}{{{\mu ^k}}}} \right\|_F^2} \right)
			      \end{array}
		      \end{equation}
		      To solve (\ref{b}), we can utilize the following system of linear equations:
		      \begin{equation}
			      \left( {2I + \Delta } \right){{\cal B}^{k + 1}} = {{L}^k_A} + {\theta _1} + {\theta _2} + {\theta _3}
		      \end{equation}
		      where
		      $\Delta  = {\cal K}_h^T{{\cal K}_h} + {\cal K}_v^T{{\cal K}_v} + {\cal K}_z^T{{\cal K}_z},{{\cal L}^k} = {\cal K} - {{\cal T}^k} - {{\cal N}^k} + \frac{{y_1^k}}{{{\mu ^k}}} + {\cal Z} + \frac{{y_2^k}}{{{\mu ^k}}},{\theta _1} = {\cal K}_h^T\left( {{\cal V}_1^k + \frac{{y_3^k}}{{{\mu ^k}}}} \right),{\theta _2} = {\cal K}_v^T\left( {{\cal V}_2^k + \frac{{y_4^k}}{{{\mu ^k}}}} \right),{\theta _3} = {\cal K}_z^T\left( {{\cal V}_3^k + \frac{{y_5^k}}{{{\mu ^k}}}} \right)$
		      ,and T is the matrix transpose.
		      By considering convolutions along two spatial directions and one temporal direction, we obtain a new Equation (\ref{newB}) as follows:
		      \begin{equation}\label{newB}
			      {{\cal B}^{k + 1}} = {F^{ - 1}}\left( {\frac{{F\left( {{{\cal L}^k} + {\theta _1} + {\theta _2} + {\theta _3}} \right)}}{{2 + {{\sum {_{i \in \left\{ {h,v,z} \right\}}F\left( {{{\cal K}_i}} \right)} }^H}F\left( {{{\cal K}_i}} \right)}}} \right)
		      \end{equation}

		\item We update the value of ${\cal T}$ using the following equation:
		      \begin{equation}\label{t}
			      \begin{split}
				      {{\cal T}^{k + 1}} &= \mathop {\arg \min }\limits_T {\lambda _s}{\left\| {\cal T} \right\|_1}\\
				      &+ \frac{{{\mu ^k}}}{2}\left\| {{\cal K} - {{\cal B}^{k + 1}} - {\cal T} - {{\cal N}^k} + \frac{{y_1^k}}{{{\mu ^k}}}} \right\|_F^2
			      \end{split}
		      \end{equation}

		      The closedform solution of (\ref{t}) can be obtained by resorting to the element-wise shrinkage operator, that is:
		      \begin{equation}
			      {{\cal T}^{k + 1}} = {\cal T}{h_{{\lambda _1}\left( {{\mu ^k}} \right)}}^{ - 1}\left( {{\cal K} - {{\cal B}^{k + 1}} - {{\cal N}^k} + \frac{{y_1^k}}{{{\mu ^k}}}} \right)
		      \end{equation}

		\item We update the values of ${\cal V}_1$, ${\cal V}_2$, and ${\cal V}_3$ using the following equations:
		      \begin{equation}
			      \label{v}
			      \left\{ \begin{array}{l}
				      \begin{split}
					      V_1^{k + 1} &= \mathop {\arg \min }\limits_{V1} {\lambda _{tv}}{\left\| {{V_1}} \right\|_1} \\
					      &+ \frac{{{\mu ^k}}}{2}\left\| {V_1^k - {{\cal K}_h}{{\cal B}^{k + 1}} + \frac{{y_3^k}}{{{\mu ^k}}}} \right\|_F^2\\
					      V_2^{k + 1} &= \mathop {\arg \min }\limits_{V2} {\lambda _{tv}}{\left\| {{V_2}} \right\|_1} \\
					      &+ \frac{{{\mu ^k}}}{2}\left\| {V_2^k - {{\cal K}_v}{{\cal B}^{k + 1}} + \frac{{y_3^k}}{{{\mu ^k}}}} \right\|_F^2\\
					      V_3^{k + 1} &= \mathop {\arg \min }\limits_{V3} \delta {\lambda _{tv}}{\left\| {{V_3}} \right\|_1} \\
					      &+ \frac{ {{\mu ^k}}}{2}\left\| {V_3^k - {{\cal K}_z}{{\cal B}^{k + 1}} + \frac{{y_3^k}}{{{\mu ^k}}}} \right\|_F^2\\
				      \end{split}
			      \end{array} \right.
		      \end{equation}
		      Here, we can also utilize the element-wise shrinkage operator to solve the above problem. The Equation is as follows:
		      \begin{equation}
			      \left\{ \begin{array}{l}
				      V_1^{k + 1} = {\cal T}{h_{{\lambda _{tv}}\left( {{\mu ^k}} \right)}}^{ - 1}\left( {{{\cal K}_h}{{\cal B}^{k + 1}} - \frac{{y_3^k}}{{{\mu ^k}}}} \right) \\
				      V_2^{k + 1} = {\cal T}{h_{{\lambda _{tv}}\left( {{\mu ^k}} \right)}}^{ - 1}\left( {{{\cal K}_v}{{\cal B}^{k + 1}} - \frac{{y_4^k}}{{{\mu ^k}}}} \right) \\
				      V_3^{k + 1} = {\cal T}{h_{\delta {\lambda _{tv}}\left( {{\mu ^k}} \right)}}^{ - 1}\left( {{{\cal K}_z}{{\cal B}^{k + 1}} - \frac{{y_5^k}}{{{\mu ^k}}}} \right)
			      \end{array} \right.
		      \end{equation}
		\item We update the value of ${{\cal N}^{k + 1}}$ using the following equation:
		      \begin{equation}\label{n}
			      \begin{split}
				      {{\cal N}^{k + 1}} &= \mathop {\arg \min }\limits_{\cal N} {\lambda _3}\left\| {\cal N} \right\|_F^2\\
				      &+ \frac{{{\mu ^k}}}{2}\left\| {{\cal K} - {{\cal B}^{k + 1}} - {{\cal T}^{k + 1}} - {\cal N} + \frac{{y_1^k}}{{{\mu ^k}}}} \right\|_F^2
			      \end{split}
		      \end{equation}

\begin{table*}
	\caption{Pictures introduction}
	\centering
	\resizebox{\linewidth}{!}
	{
		\begin{tabular}{c|c|c|c|c|c}
			\toprule
			\toprule
			Sequence & Frames & Image Size     & Average SCR & Target Descriptions                             & Background Descriptions \\
			\hline
			1        & 120    & $256\times256$ & 5.45        & Far range, single target                        & ground                  \\ 
			2        & 120    & $256\times256$ & 5.11        & Near to far, single target                      & ground                  \\ 
			3        & 120    & $256\times256$ & 6.33        & Near to far, single target                      & ground                  \\ 
			4        & 120    & $256\times256$ & 6.07        & Far to near, single target                      & ground                  \\ 
			5        & 120    & $256\times256$ & 5.20        & Far distance, single target                     & ground-sky boundary     \\ 
			6        & 120    & $256\times256$ & 1.98        & Target near to far, single target, faint target & ground                  \\ 
			\bottomrule
			\toprule
		\end{tabular}
	}
	\label{Datasets_INTRODUCTION2}
\end{table*}
\begin{table*}
	\caption{Parameter setting}
	\centering
	\resizebox{\linewidth}{!}
	{
		\begin{tabular}{c|c|c}
			\toprule
			\toprule
			Methods                                                            & Acronyms      & Parameter settings                                                                                             \\
			\hline
			Total Variation Regularization and Principal Component Pursuit \cite{tvpcp}      & TV-PCP       & $\lambda_{1}=0.005$, $\lambda_{2}=\frac{1}{\sqrt{max(M,N)}}$, $\beta=0.025$, $\gamma=1.5$                      \\
			Partial Sum of the Tensor Nuclear Norm \cite{PSTNN}                             & PSTNN         & Sliding step: 40, $\lambda=\frac{0.6}{\sqrt{max(n1,n2)}}, pathch size:30\times30$                              \\
			Multiple Subspace Learning and Spatial-temporal Patch-Tensor Model \cite{MSLSTIPT} & MSLSTIPT      & $L=6$,$P=0.8$,$\lambda=\frac{1}{\sqrt{n_{3}\times max(n_{1},n_{2})}}$,$patch size:30\times 30$                 \\
			Nonconvex tensor fibered rank approximation \cite{NTFRA}                        & NTFRA         & $Sliding step:40$,$patch size:40\times40$,$\lambda=\frac{1}{\sqrt{n_{3}\times max(n_{1},n_{2})}}$,$\beta=0.01$ \\
			Non-Convex Tensor Low-Rank Approximation \cite{ASSTV}                           & ASSTV-NTLA    & $L=3$, $H=6$, $\lambda_{tv}=0.5$, $\lambda_{s}=\frac{H}{\sqrt{max(M, N)\times L}}$, $\lambda_{3}=100$           \\
			Tensor $L_{2,1}$ Norm Minimization via Tensor QR Decomposition     & ASSTV-TLNMTQR & $r=180$, $L=3$, $H=6$, $\lambda_{tv}=0.5$, $\lambda_{s}=\frac{H}{\sqrt{max(M, N)\times L}}$, $\lambda_{3}=100$  \\
			\bottomrule
			\toprule
		\end{tabular}
	}
	\label{PARAMETER_SETTING}
\end{table*}

		      The solution of (\ref{n}) is as follows:
		      \begin{equation}
			      {{\cal N}^{k + 1}} = \frac{{{\mu ^k}\left( {{\cal K} - {{\cal B}^{k + 1}} - {{\cal T}^{k + 1}}} \right) + y_1^k}}{{{\mu ^k} + 2{\lambda _3}}}
		      \end{equation}
		\item Updating multipliers y1; y2; y3; y4; y5 with other variables being fixed:
		      \begin{equation}
			      \label{y}
			      \left\{ \begin{array}{l}
				      y_1^{k + 1} = y_1^k + {\mu ^k}\left( {{\cal K} - {{\cal B}^{k + 1}} - {{\cal T}^{k + 1}} - {{\cal N}^{k + 1}}} \right) \\
				      y_2^{k + 1} = y_2^k + {\mu ^k}\left( {{{\cal Z}^{k + 1}} - {{\cal B}^{k + 1}}} \right)                                 \\
				      y_3^{k + 1} = y_3^k + {\mu ^k}\left( {V_1^{k + 1} - {{\cal K}_h}{{\cal B}^{k + 1}}} \right)                            \\
				      y_4^{k + 1} = y_4^k + {\mu ^k}\left( {V_2^{k + 1} - {{\cal K}_v}{{\cal B}^{k + 1}}} \right)                            \\
				      y_5^{k + 1} = y_5^k + {\mu ^k}\left( {V_3^{k + 1} - {{\cal K}_z}{{\cal B}^{k + 1}}} \right)
			      \end{array} \right.
		      \end{equation}
		\item Updating ${\mu ^{k + 1}}$ by ${\mu ^{k + 1}} = \min (\rho {\mu ^k},{\mu _{\max }})$.

		      \begin{algorithm}[H]
                    \small
			      \caption{: ASSTV-TLNMTQR algorithm}
			      \label{algorithm2}
			      \begin{algorithmic}[1]
				      \REQUIRE infrared image sequence ${k_1},{k_2}, \cdots ,{k_P} \in \mathbb{R}  {^{{n_1} \times {n_2}}}$,\quad number of frames $L$,\quad parameters ${\lambda _s},\quad {\lambda _{tv}},\quad {\lambda _3},\quad \mu \textgreater 0 $;
                        \STATE Initialize: Transform the image sequence into the original tensor ${\cal K}$,\quad${\cal B}^{0}$=${\cal T}^{0}$=${\cal N}^{0}$=${\cal V}_{i}^{0}$=0,\quad$i=1,2,3$,\quad${y}_{i}^0$=0,\quad${\mu}_{0}=0.005$,\quad${\mu}_{max}=1e7$,\quad$k=0$,\quad$i=1,\dots,5$,\quad$\rho=1,5$,\quad$\xi=1e-6$.
				      \WHILE{Not converged}
				      \STATE Update ${{\cal Z}^{k + 1}}$ by \ref{z}
				      \STATE Update ${{\cal B}^{k + 1}}$ by \ref{b}
				      \STATE Update ${{\cal T}^{k + 1}}$ by \ref{t}
				      \STATE Update $V_1^{k + 1}$,$V_2^{k + 1}$,$V_3^{k + 1}$ by \ref{v}
				      \STATE Update ${{\cal N}^{k + 1}}$ by \ref{n}
				      \STATE Update multipliers $y_i^{k + 1}$, $i=1,2,\cdots,5 $ by \ref{y}
				      \STATE Updating ${\mu ^{k + 1}}$ by ${\mu ^{k + 1}} = \min (\rho {\mu ^k},{\mu _{\max }})$
				      \STATE  Check the convergence conditions
				      $\frac{{\left\| {{\cal K} - {{\cal B}^{k + 1}} - {{\cal T}^{k + 1}} - {\cal N} + \frac{{y_1^k}}{{{\mu ^k}}}} \right\|_F^2}}{{\left\| {\cal K} \right\|_F^2}} \leq \xi $
				      \STATE Update $k=k+1$
				      \ENDWHILE
				      \ENSURE  ${{\cal B}^{k + 1}}$,${{\cal T}^{k + 1}}$,${{\cal N}^{k + 1}}$
			      \end{algorithmic}
		      \end{algorithm}
	\end{enumerate}

\section{NUMERICAL TESTS AND APPLICATIONS}
	The main purpose of our algorithm is to decompose the tensor $ \cal K$ into three tensors: the background tensor $\cal B$, the target tensor $\cal T$, and the noise tensor $\cal N$. We will explore the effects of different parameter adjustments and ultimately identify a relatively suitable parameter configuration for the most effective extraction of infrared small target detection. Finally, we will compare our algorithm with the latest algorithms.

\begin{figure*}[!htb]
	\centering
	\begin{minipage}{0.3\textwidth}
		\centering
		\caption{The parameter analysis of r=10}
		\includegraphics[width=0.9\linewidth]{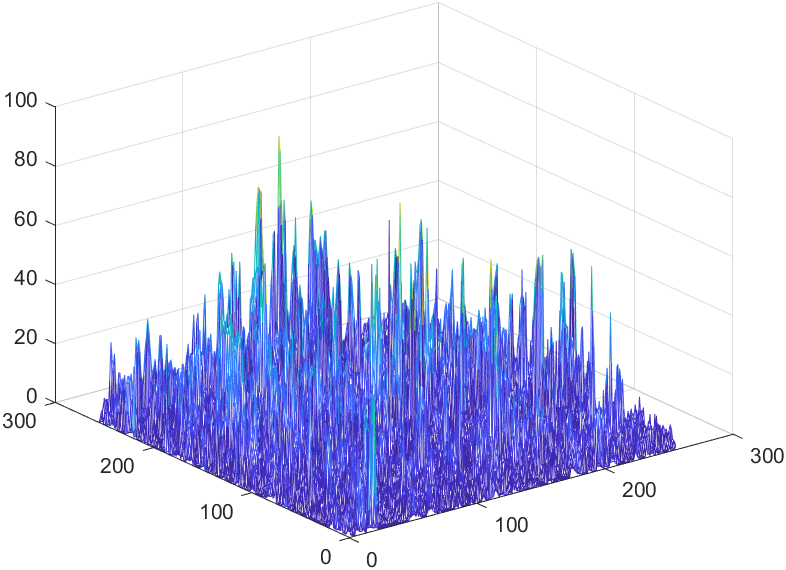}
		\label{r1}
	\end{minipage}
	\begin{minipage}{0.3\textwidth}
		\centering
		\caption{The parameter analysis of r=50}
		\includegraphics[width=0.9\linewidth]{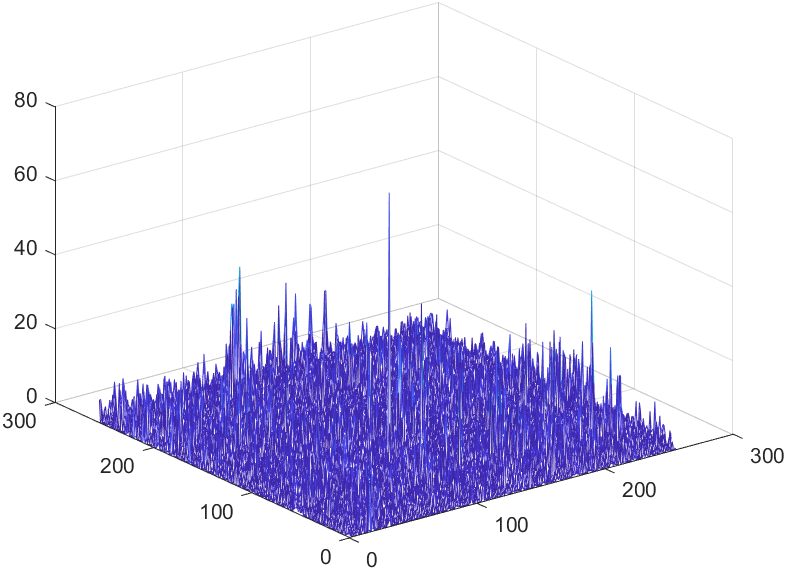}
		\label{r2}
	\end{minipage}
	\begin{minipage}{0.3\textwidth}
		\centering
		\caption{The parameter analysis of r=90}
		\includegraphics[width=0.9\linewidth]{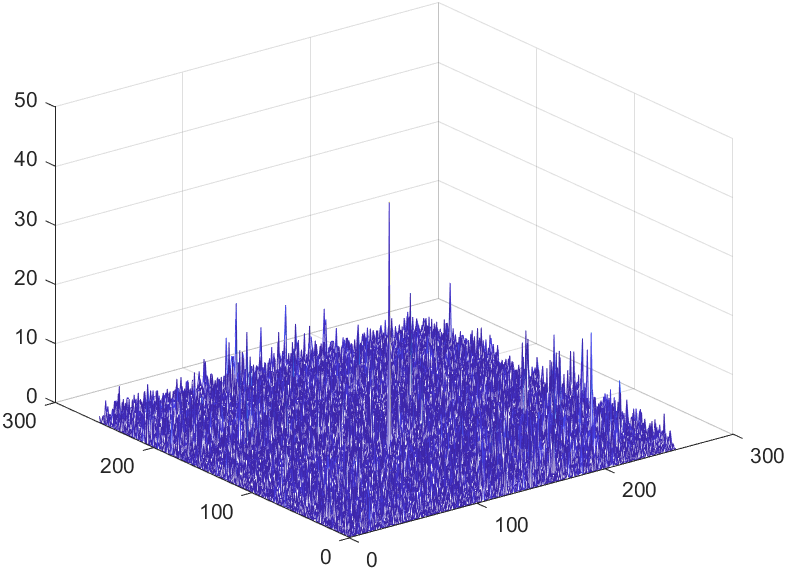}
		\label{r3}
	\end{minipage}

	\begin{minipage}{0.3\textwidth}
		\centering
		\caption{The parameter analysis of r=130}
		\includegraphics[width=0.9\linewidth]{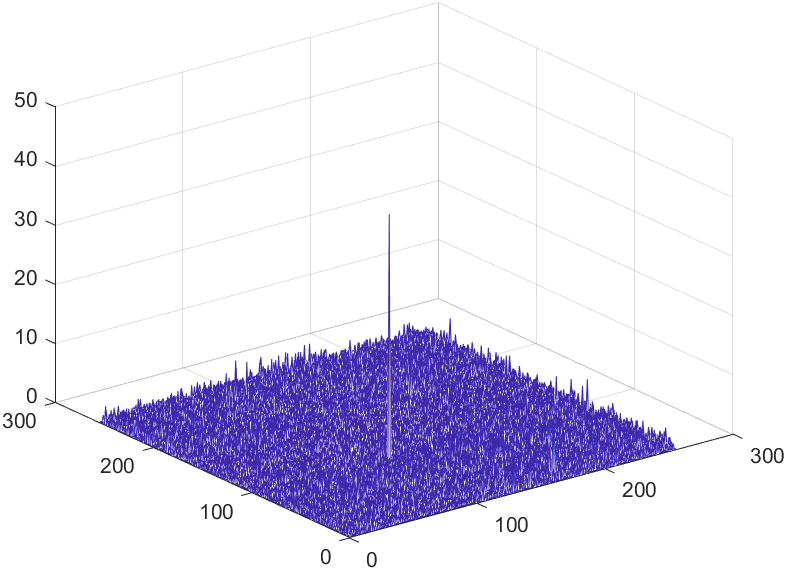}
		\label{r4}
	\end{minipage}
	\begin{minipage}{0.3\textwidth}
		\centering
		\caption{The parameter analysis of r=170}
		\includegraphics[width=0.9\linewidth]{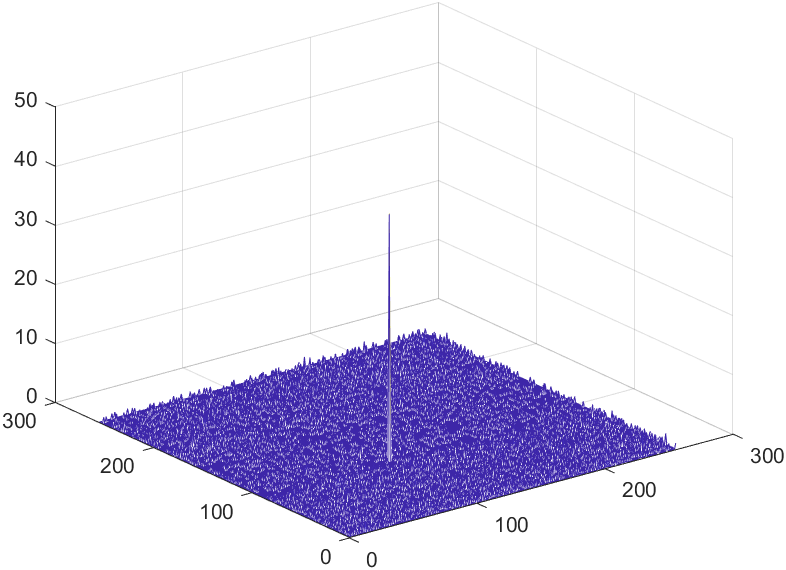}
		\label{r5}
	\end{minipage}
	\begin{minipage}{0.3\textwidth}
		\centering
		\caption{The parameter analysis of r=210}
		\includegraphics[width=0.9\linewidth]{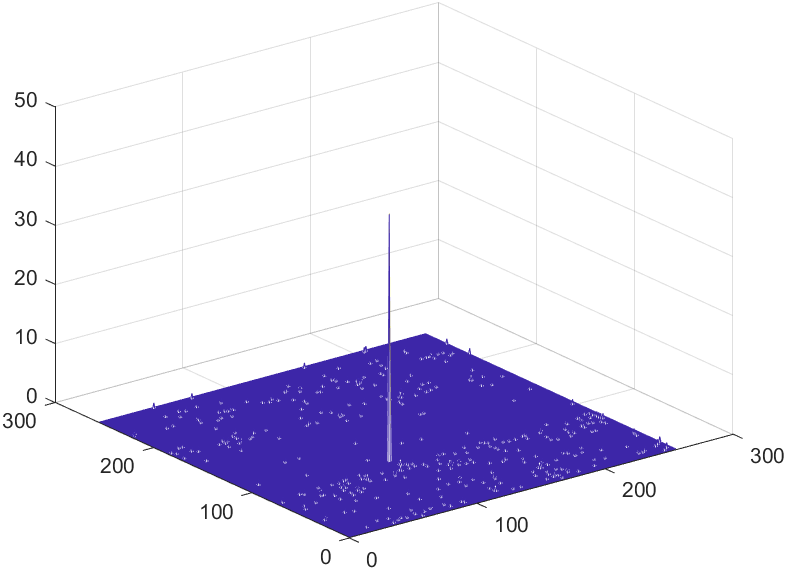}
		\label{r6}
	\end{minipage}
\end{figure*}
\begin{figure*}[!htb]
	\centering
	\begin{minipage}{0.3\textwidth}
		\centering
		\caption{Parameter analysis of $L$ in Sequence1}
		\includegraphics[width=0.9\linewidth]{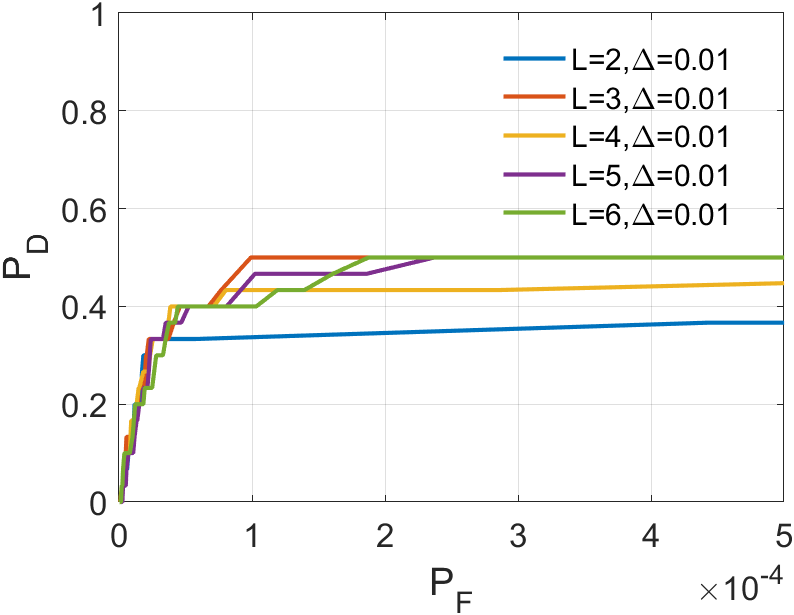}
		\label{sl1}
	\end{minipage}
	\begin{minipage}{0.3\textwidth}
		\centering
		\caption{Parameter analysis of $L$ in Sequence2}
		\includegraphics[width=0.9\linewidth]{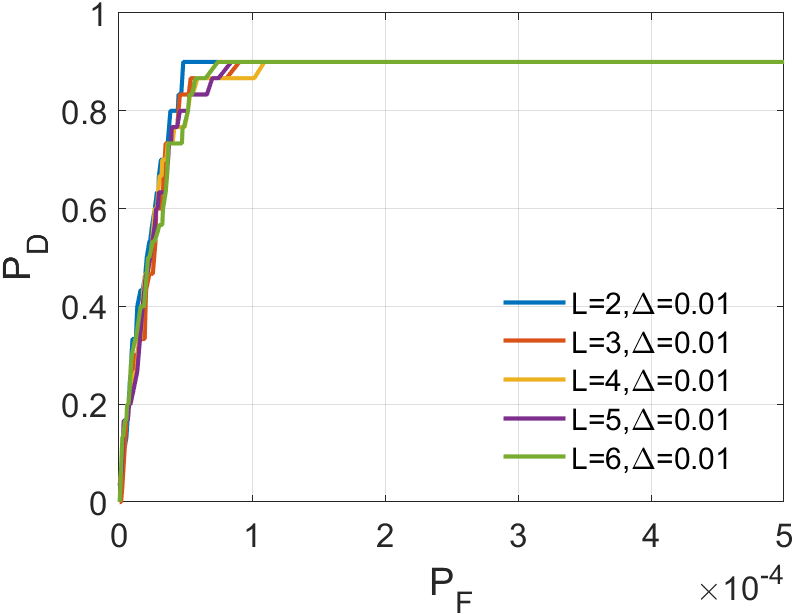}
		\label{sl2}
	\end{minipage}
	\begin{minipage}{0.3\textwidth}
		\centering
		\caption{Parameter analysis of $L$ in Sequence3}
		\includegraphics[width=0.9\linewidth]{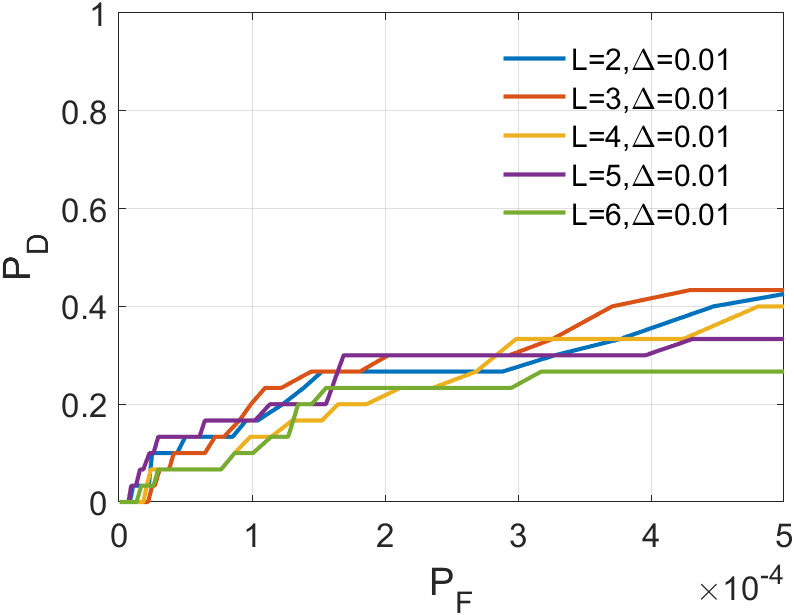}
		\label{sl3}
	\end{minipage}

	\begin{minipage}{0.3\textwidth}
		\centering
		\caption{Parameter analysis of $L$ in Sequence4}
		\includegraphics[width=0.9\linewidth]{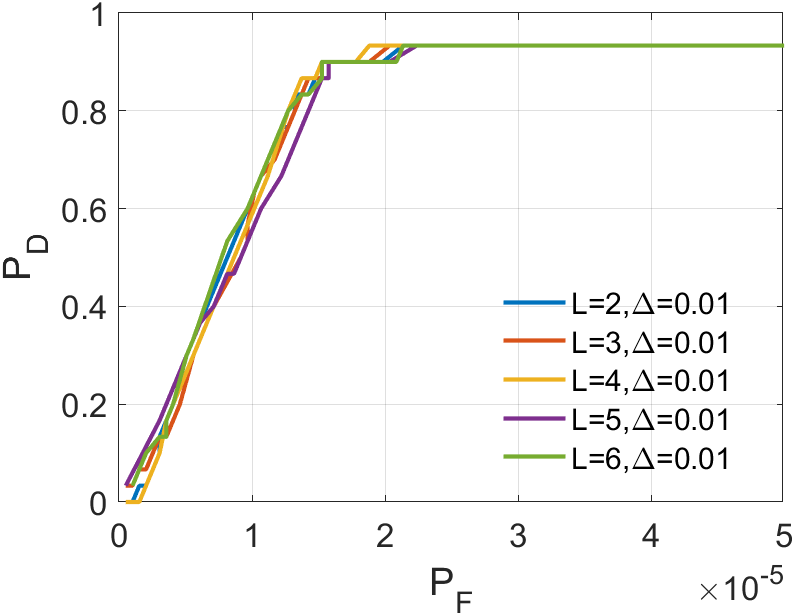}
		\label{sl4}
	\end{minipage}
	\begin{minipage}{0.3\textwidth}
		\centering
		\caption{Parameter analysis of $L$ in Sequence5}
		\includegraphics[width=0.9\linewidth]{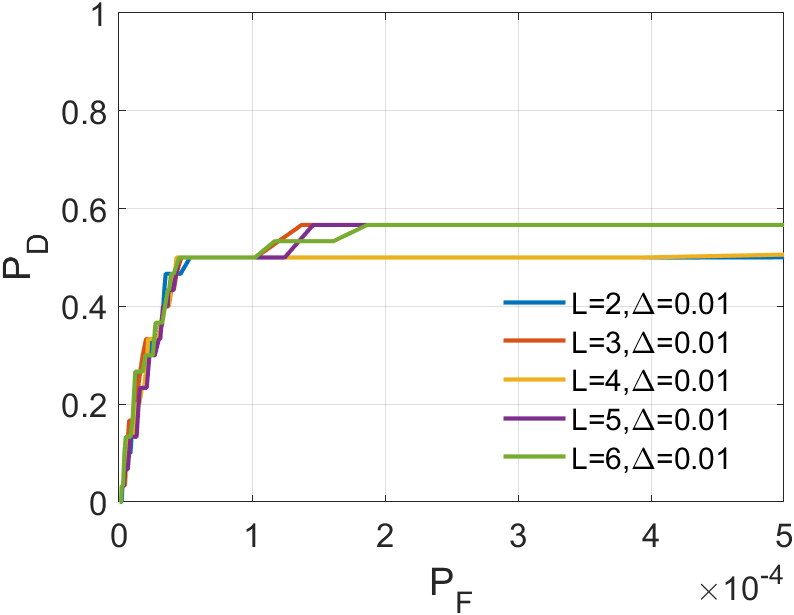}
		\label{sl5}
	\end{minipage}
	\begin{minipage}{0.3\textwidth}
		\centering
		\caption{Parameter analysis of $L$ in Sequence6}
		\includegraphics[width=0.9\linewidth]{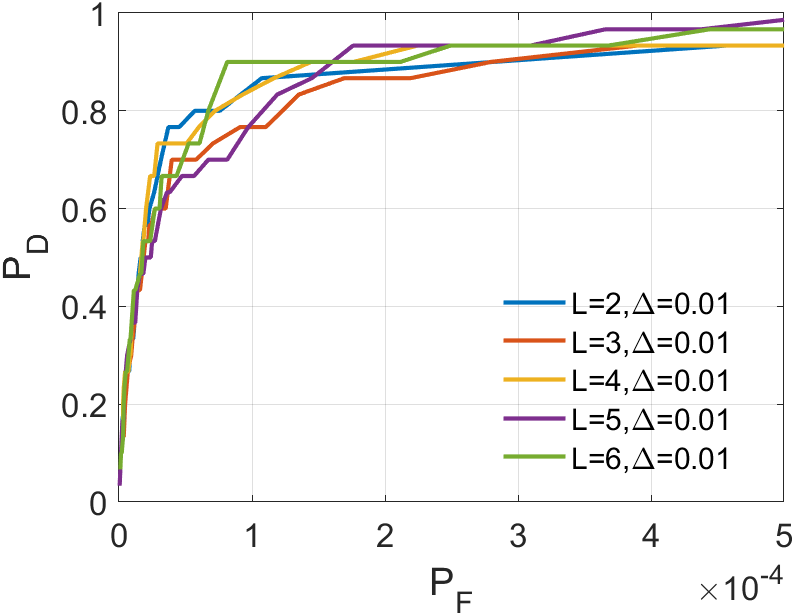}
		\label{sl6}
	\end{minipage}
\end{figure*}

\subsection{Infrared small target detection}

This section aims to demonstrate the robustness of our algorithm through a series of experiments related to infrared small target detection. We employ 3D Receiver Operating Characteristic Curve (ROC) \cite{ROC}  to assess the performance of our algorithm in separating tensor $\mathcal{T}$. Additionally, we have selected five comparative methods to highlight the superiority of our approach.

\subsubsection{Evaluation metrics and baseline methods}

In the experiments related to small infrared targets, we'll use 3D ROC to assess the algorithm's capability. The 3D ROC are based on 2D ROC \cite{ROC2, 2D_ROC}. Initially, we plot the $(PD, PF)$ curve, followed by separate curves for $(PD, \tau)$ and $(PF, \tau)$, and finally combine them to form the 3D ROC.

The horizontal axis of 2D ROC represents the false alarm rate \(F_a\), and the vertical axis represents the detection probability \(P_d\). They are defined as follows:
\begin{equation}
	\label{pd}
	{P_d} = \frac{{{\rm{number\ of\ true\ detections}}}}{{{\rm{number\ of\  actual\  targets}}}}
\end{equation}
\begin{equation}
	\label{fa}
	{F_a} = \frac{{{\rm{number\  of\  false \ detections}}}}{{{\rm{number\  of\  image\  pixels}}}}
\end{equation}
The above two indicators range between 0 and 1.

\begin{figure*}[!htb]
	\centering
	\begin{minipage}{0.3\textwidth}
		\centering
		\caption{Parameter analysis of $H$ in Sequence1}
		\includegraphics[width=0.9\linewidth]{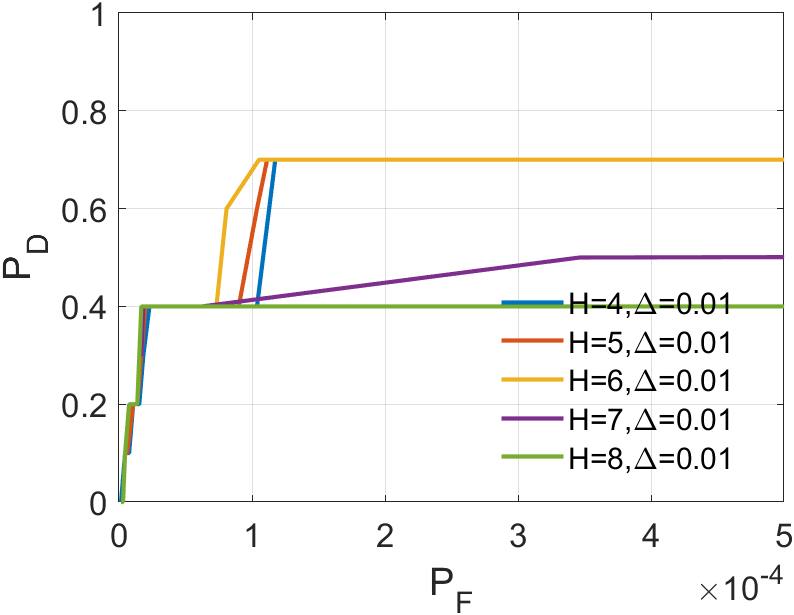}
		\label{sh1}
	\end{minipage}
	\begin{minipage}{0.3\textwidth}
		\centering
		\caption{Parameter analysis of $H$ in Sequence2}
		\includegraphics[width=0.9\linewidth]{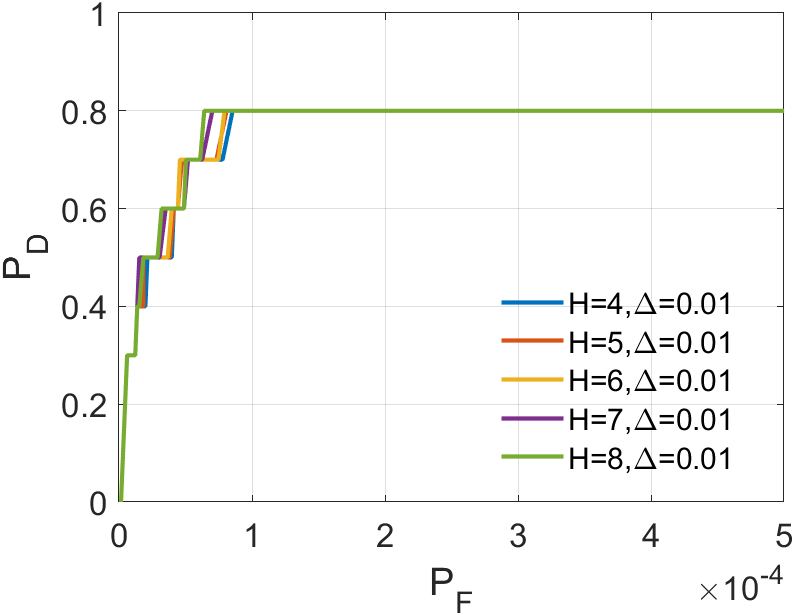}
		\label{sh2}
	\end{minipage}
	\begin{minipage}{0.3\textwidth}
		\centering
		\caption{Parameter analysis of $H$ in Sequence3}
		\includegraphics[width=0.9\linewidth]{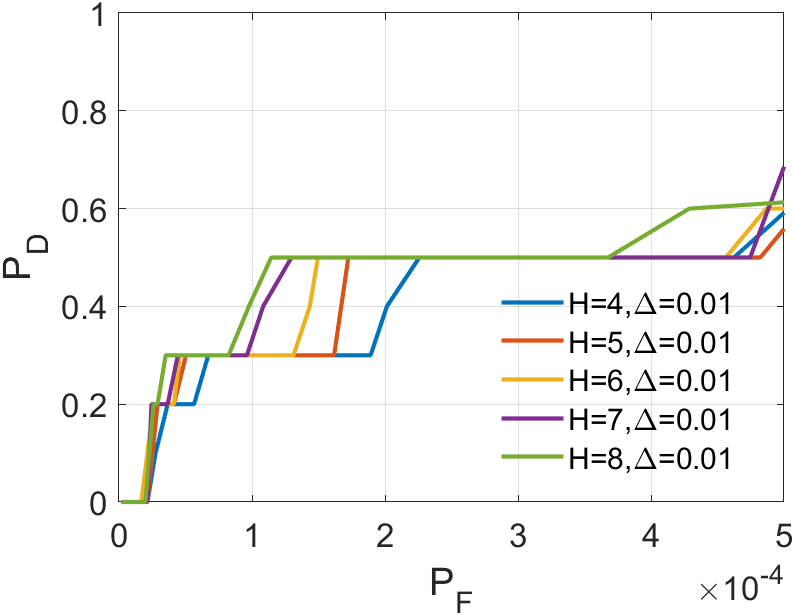}
		\label{sh3}
	\end{minipage}

	\begin{minipage}{0.3\textwidth}
		\centering
		\caption{Parameter analysis of $H$ in Sequence4}
		\includegraphics[width=0.9\linewidth]{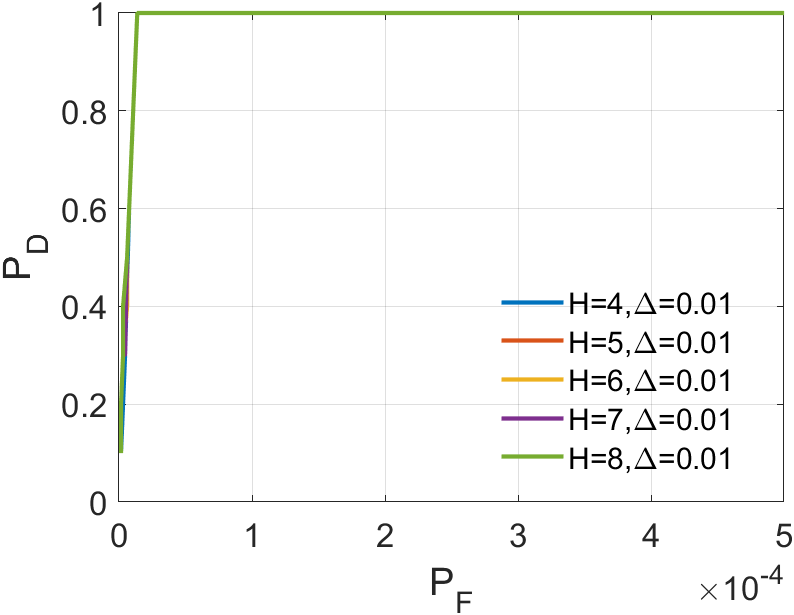}
		\label{sh4}
	\end{minipage}
	\begin{minipage}{0.3\textwidth}
		\centering
		\caption{Parameter analysis of $H$ in Sequence5}
		\includegraphics[width=0.9\linewidth]{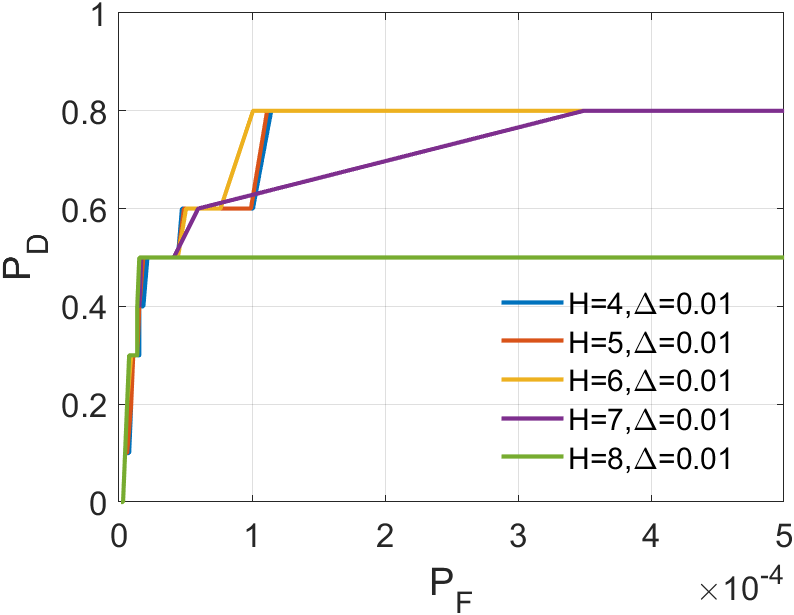}
		\label{sh5}
	\end{minipage}
	\begin{minipage}{0.3\textwidth}
		\centering
		\caption{Parameter analysis of $H$ in Sequence6}
		\includegraphics[width=0.9\linewidth]{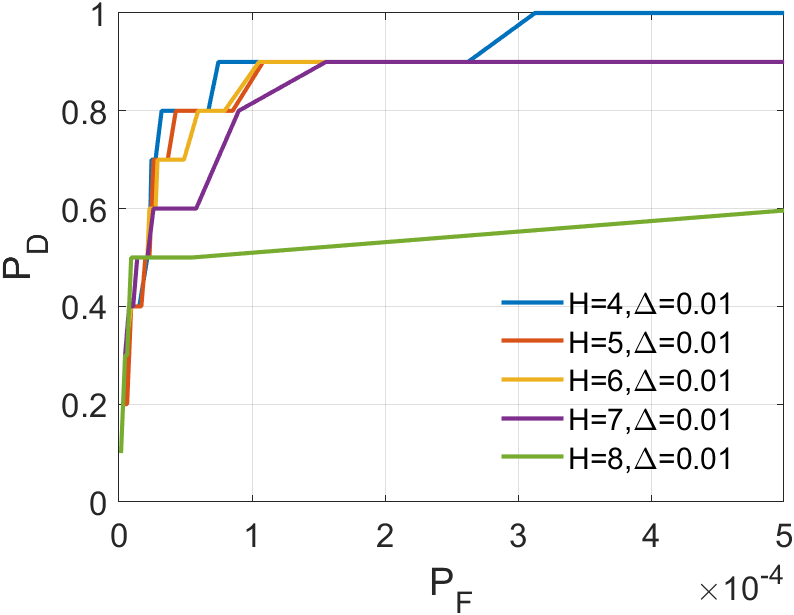}
		\label{sh6}
	\end{minipage}
\end{figure*}
\begin{table*}
	\caption{Quantitative comparasion of different methods on sequences 1-6}
	\centering
	\resizebox{\linewidth}{!}
	{
		\begin{tabular}{c|c|c|c|c|c|c|c}
			\toprule
			\toprule
			Sequence & \diagbox{Metrics}{Method} & ASSTV-TLNMTQR(ours) & ASSTV-NTLA & NTFRA  & MSLSTIPT & PSTNN & TVPCP \\

			\hline
			\multirow{3}{*}{Sequence 1}
			         & $AUC(PF,PD)$                   & 0.844               & 0.849      & 0.881  & 0.701    & 0.682 & 0.848 \\
			         & $AUC(PF,\tau)$                 & 0.005               & 0.005      & 0.054  & 0.013    & 0.021 & 0.006 \\
			         & $RUN TIME(s)$               & 4.135               & 5.304      & 25.688 & 0.161    & 3.031 & 0.490 \\
			\hline
			\multirow{3}{*}{Sequence 2}
			         & $AUC(PF,PD)$                    & 0.895               & 0.899      & 0.965  & 0.833    & 0.998 & 0.950 \\
			         & $AUC(PF,\tau)$                  & 0.005               & 0.005      & 0.054  & 0.011    & 0.009 & 0.006 \\
			         & $RUN TIME(s) $              & 4.337               & 5.558      & 27.126 & 0.219    & 3.112 & 0.442 \\
			\hline
			\multirow{3}{*}{Sequence 3}
			         & $AUC(PF,PD)$                    & 0.895               & 0.899      & 0.960  & 0.686    & 0.946 & 0.899 \\
			         & $AUC(PF,\tau)$                 & 0.005               & 0.005      & 0.068  & 0.008    & 0.015 & 0.006 \\
			         & $RUN TIME(s) $                 & 4.476               & 5.806      & 26.409 & 0.251    & 2.948 & 0.369 \\
			\hline
			\multirow{3}{*}{Sequence 4}
			         & $AUC(PF,PD)$                    & 0.998               & 0.998      & 0.951  & 0.947    & 0.972 & 0.981 \\
			         & $AUC(PF,\tau)$                & 0.006               & 0.005      & 0.59   & 0.010    & 0.023 & 0.007 \\
			         & $RUN TIME(s) $              & 4.061               & 5.188      & 21.649 & 0.213    & 2.732 & 0.443 \\
			\hline
			\multirow{3}{*}{Sequence 5}
			         & $AUC(PF,PD)$                    & 0.896               & 0.8999     & 0.951  & 0.705    & 0.781 & 0.899 \\
			         & $AUC(PF,\tau)$                & 0.005               & 0.005      & 0.054  & 0.014    & 0.022 & 0.008 \\
			         & $RUN TIME(s)$               & 4.136               & 5.133      & 24.497 & 0.159    & 2.817 & 0.517 \\
			\hline
			\multirow{3}{*}{Sequence 6}
			         & $AUC(PF,PD)$                    & 0.947               & 0.949      & 0.995  & 0.937    & 0.645 & 0.750 \\
			         & $AUC(PF,\tau)$                  & 0.006               & 0.005      & 0.066  & 0.013    & 0.009 & 0.006 \\
			         & $RUN TIME(s) $              & 4.367               & 5.203      & 25.966 & 0.195    & 2.869 & 0.492 \\
			\hline
			\bottomrule
		\end{tabular}
	}
	\label{comp}
\end{table*}

In all ROC, we assess the algorithm's performance by comparing the Area Under Curve (AUC). In 2D ROC, a bigger AUC generally signifies better algorithm performance. However, when the algorithm's false alarm rate is high, indicating a higher detection of irrelevant points, it can create a false impression of a large area under the $(PD, PF)$ curve. Hence, the 2D ROC may not fully reflect the algorithm's performance, prompting our choice of the 3D ROC. Within the 3D ROC, we've incorporated the relationships between PD, PF, and the threshold $\tau$. The curve $(PD, \tau)$ has the same implication as $(PD, PF)$, a larger area signifies better performance. Therefore, for subsequent experiments, we will solely showcase the $(PD, PF)$ curve. Meanwhile, the $(PD, \tau)$ curve reflects the algorithm's background suppression capability, specifically the quantity of irrelevant points. Generally, a smaller area under the $(PD, \tau)$ curve implies a stronger background suppression ability and better algorithm performance.

\subsubsection{Parameter setting and datasets}

Above, we provide detailed parameter settings for our method and the comparative method, as shown in TABLE \ref{PARAMETER_SETTING}. We selected six different sets of images from the dataset \cite{istd_data}, composing a time sequence of 120 frames. For specific details of the images, refer to TABLE \ref{Datasets_INTRODUCTION2}.

\subsection{Parameter analysis}
\subsubsection{$r$ of tensor $\cal D$}
In our algorithm, the tensor $\mathcal{D} \in \mathbb{C}^{r \times r \times n_3}$ serves as an approximation of tensor $\mathcal{Z}$ in Equation (\ref{z}), making the size $r$ of tensor $\mathcal{D}$ the most critical parameter. The variation in $\mathcal{D}$'s size significantly impacts our algorithm's performance in detecting weak infrared targets. Refer to Fig. \ref{r1} to Fig. \ref{r6}. We utilized Sequence 6 to run our algorithm. By plotting mesh grids of different $r$ values on the same frame, a clear trend emerges: as $r$ gradually increases, the detection performance for infrared small targets improves. Additionally, in the mesh grids for $r=170$ and $r=210$, the differences in detection results become negligible. To pursue faster processing, we have chosen $r=180$ as the parameter for our subsequent experiments.

\subsubsection{Number of frames $L$}
In a sequence, the number of frames inputted each time significantly affects the $n3$ dimension of the tensor we construct. The $n3$ dimension primarily encapsulates temporal information, greatly influencing the accuracy of ASSTV calculations. Hence, we need to test the frame quantity to identify the optimal number of frames. In our algorithm, the frame quantity is denoted by the variable $L$. We vary $L$ from 2 to 6 with an increment of 1.
From the experimental results in \ref{sl1} to \ref{sl6}, we observe that in Sequence 1, Sequence 3, and Sequence 5, the best results are achieved when $L=3$, with minimal differences in the remaining sequences. After seeing the false alarm rates of the remaining sequences, we found that $L=3$ is the lowest. Therefore, for subsequent experiments, we will adopt $L=3$.

\begin{figure*}[!htb]
	\caption{The 3D ROC of Sequence 1}
	\centering
	\label{sc1}
	\begin{minipage}{0.3\textwidth}
		\centering
		\includegraphics[width=0.9\linewidth]{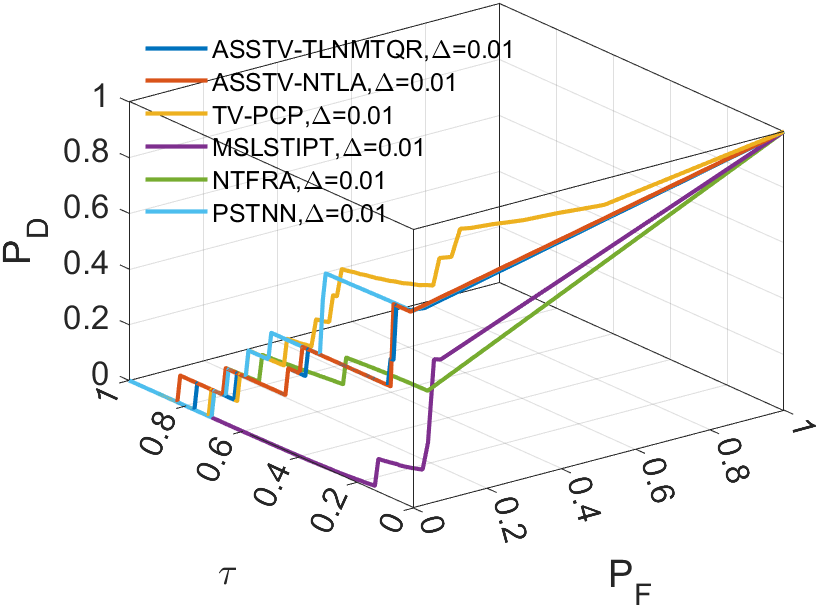}
		\label{sc1-1}
	\end{minipage}
	\begin{minipage}{0.3\textwidth}
		\centering
		\includegraphics[width=0.9\linewidth]{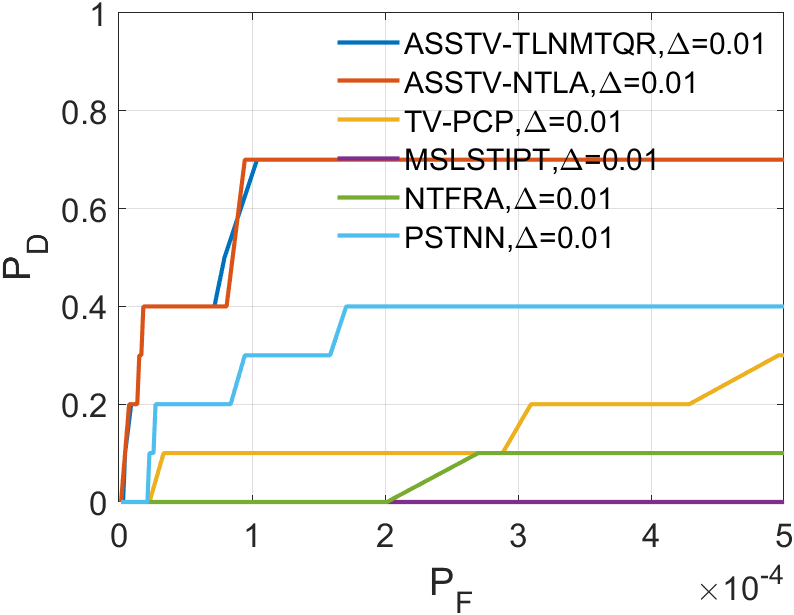}
		\label{sc1-2}
	\end{minipage}
	\begin{minipage}{0.3\textwidth}
		\centering
		\includegraphics[width=0.9\linewidth]{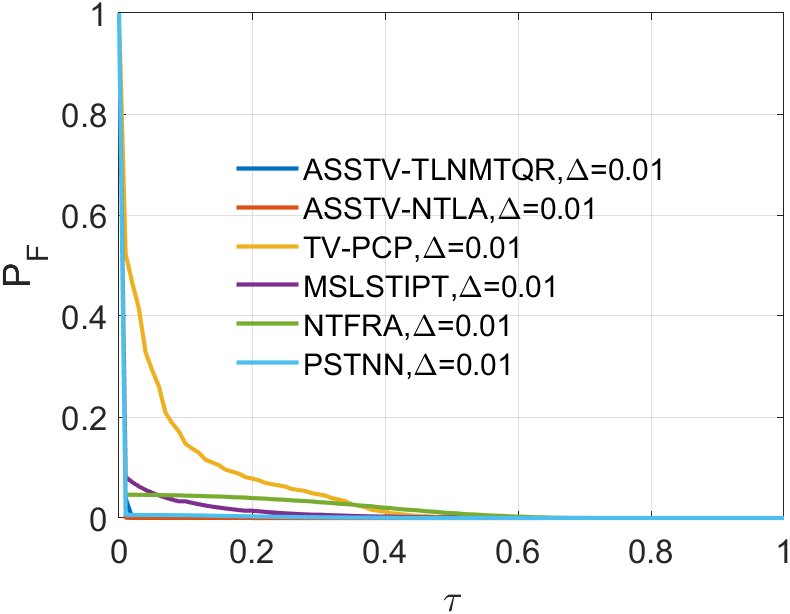}
		\label{sc1-3}
	\end{minipage}
\end{figure*}
\begin{figure*}[!htb]
	\caption{The 3D ROC of Sequence 2}
	\centering
	\label{sc2}
	\begin{minipage}{0.3\textwidth}
		\centering
		\includegraphics[width=0.9\linewidth]{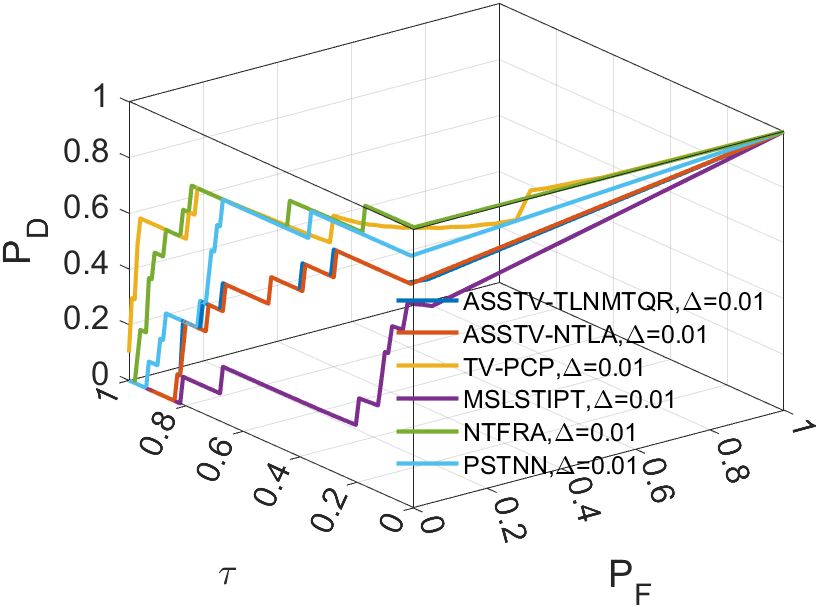}
		\label{sc2-1}
	\end{minipage}
	\begin{minipage}{0.3\textwidth}
		\centering
		\includegraphics[width=0.9\linewidth]{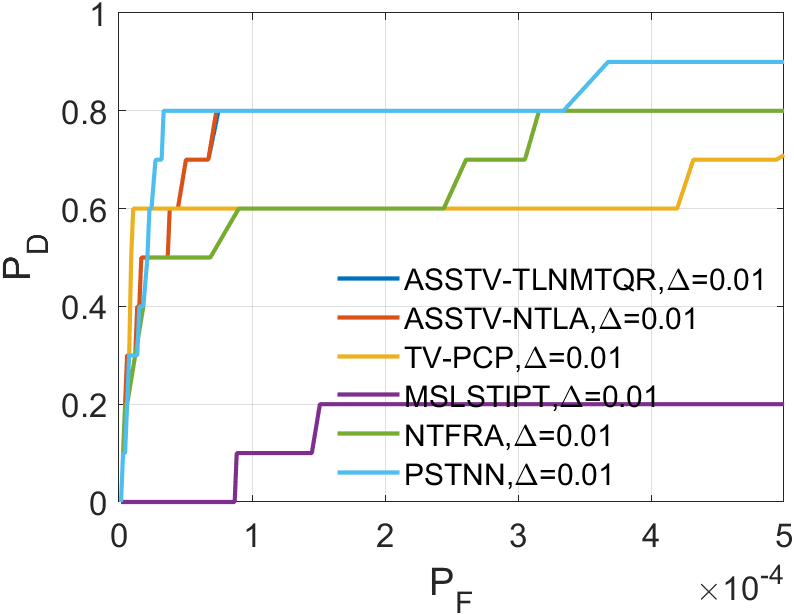}
		\label{sc2-2}
	\end{minipage}
	\begin{minipage}{0.3\textwidth}
		\centering
		\includegraphics[width=0.9\linewidth]{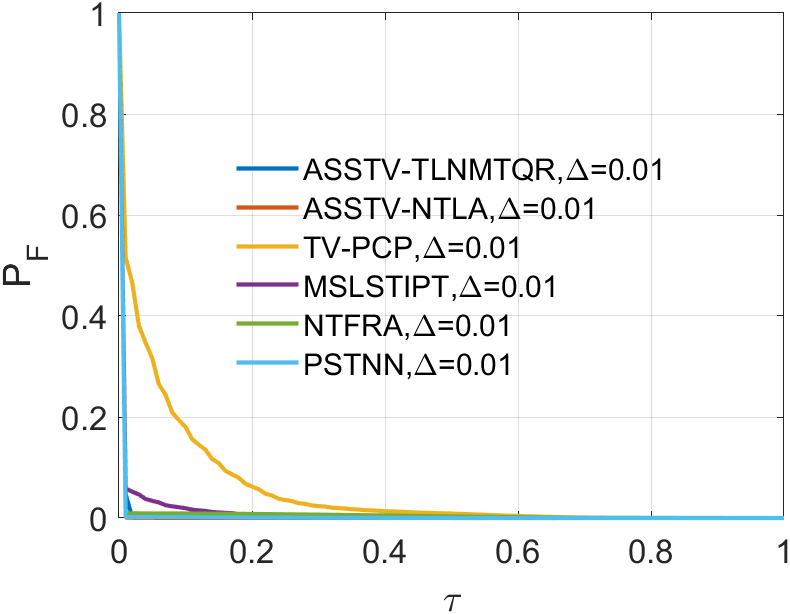}
		\label{sc2-3}
	\end{minipage}
\end{figure*}
\begin{figure*}[!htb]
	\caption{The 3D ROC of Sequence 3}
	\centering
	\label{sc3}
	\begin{minipage}{0.3\textwidth}
		\centering
		\includegraphics[width=0.9\linewidth]{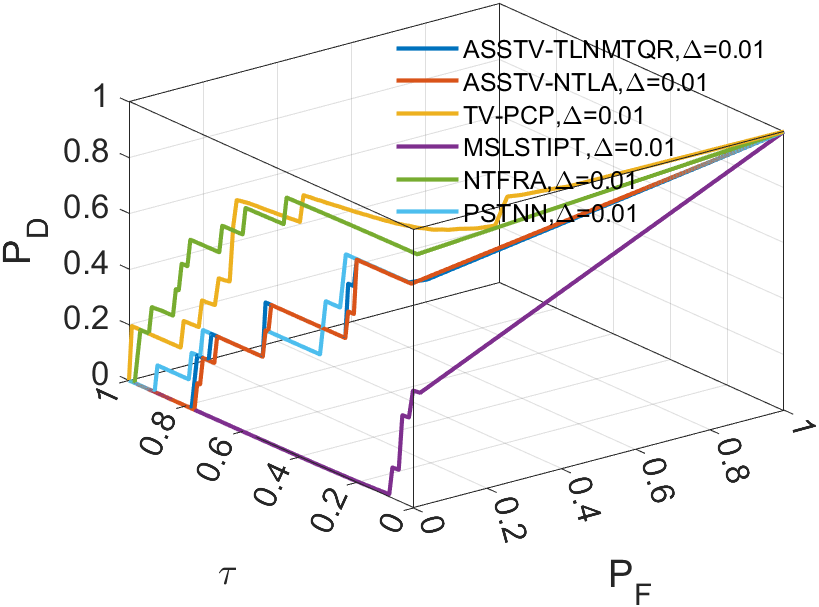}
		\label{sc3-1}
	\end{minipage}
	\begin{minipage}{0.3\textwidth}
		\centering
		\includegraphics[width=0.9\linewidth]{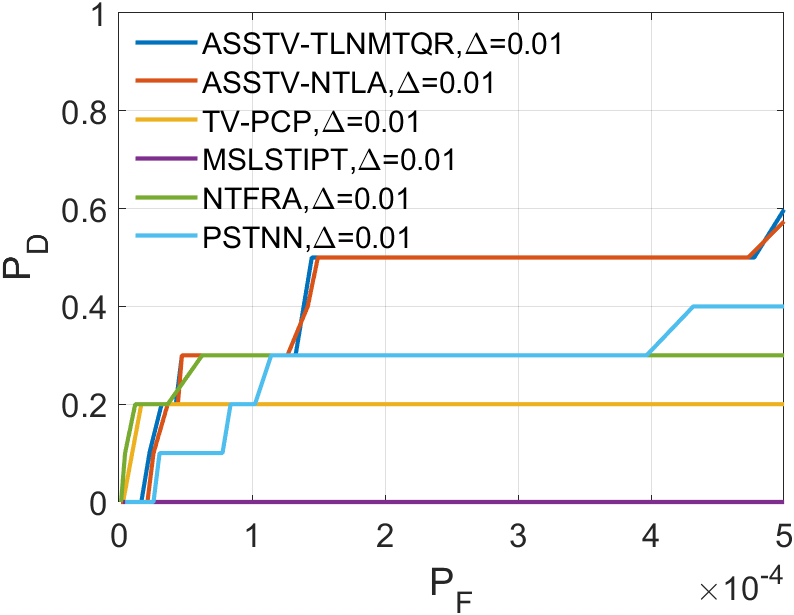}
		\label{sc3-2}
	\end{minipage}
	\begin{minipage}{0.3\textwidth}
		\centering
		\includegraphics[width=0.9\linewidth]{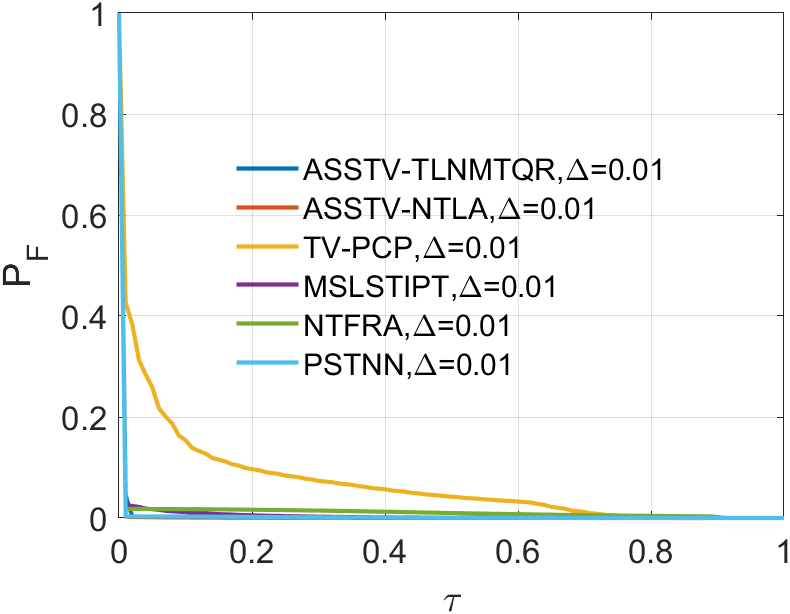}
		\label{sc3-3}
	\end{minipage}
\end{figure*}

\subsubsection{Tunning parameter $H$}
$H$ plays a crucial role in optimizing the ASTTV-TLNMTQR model. We vary $H$ from 2 to 10 with increments of 2. From the results in \ref{sh1} to \ref{sh6}, it's evident that $H$ has different effects across various backgrounds. In Sequence 1 and Sequence 5, $H=6$ emerges as the best choice, maintaining a top-three position in other sequences as well. Although $H=8$ performs best in Sequence 2 and Sequence 3, but its false alarm rate is high, so its result in Sequence 2 and Sequence 3 means nothing. Hence, we consider $H=6$ as the most stable parameter, and we will continue using this value in the subsequent experiments.

\subsection{Comparison to state-of art methods}
In our comparative experiments, we systematically evaluated the performance of our ASSTV-NTLA algorithm against five other algorithms, all of which were tested on the six sequences outlined in TABLE \ref{Datasets_INTRODUCTION2}. The specific parameter settings for each algorithm can be found in \ref{PARAMETER_SETTING}.

Upon analyzing the experimental results depicted in Fig. \ref{sc1} to Fig. \ref{sc2}, we observed that certain algorithms demonstrated detection rates either higher or nearly comparable to ASSTV-NTLA and ASSTV-TLNMTQR (ours) in specific sequences. However, a closer examination of TABLE \ref{comp} reveals that these seemingly impressive detection rates are accompanied by high false alarm rates. This critical insight underscores the superior performance of ASSTV-NTLA and ASSTV-TLNMTQR (ours), as they consistently achieve the highest detection rates while maintaining lower false alarm rates. Furthermore, the comprehensive analysis presented in TABLE \ref{comp} highlights an additional advantage of our ASSTV-TLNMTQR algorithm, which consistently performs at a speed approximately $25\%$ faster than ASSTV-NTLA across all six sequences. Notably, this accelerated speed is achieved without compromising accuracy, as evidenced by the mere $0.2\%$ difference in average accuracy compared to the ASSTV-NTLA method.

In summary, our ASSTV-TLNMTQR algorithm is the optimal choice. It has the fastest speed while ensuring higher detection rate and lower false alarm rate.

\begin{figure*}[!htb]
	\caption{The 3D ROC of Sequence 4}
	\centering
	\label{sc4}
	\begin{minipage}{0.3\textwidth}
		\centering
		\includegraphics[width=0.9\linewidth]{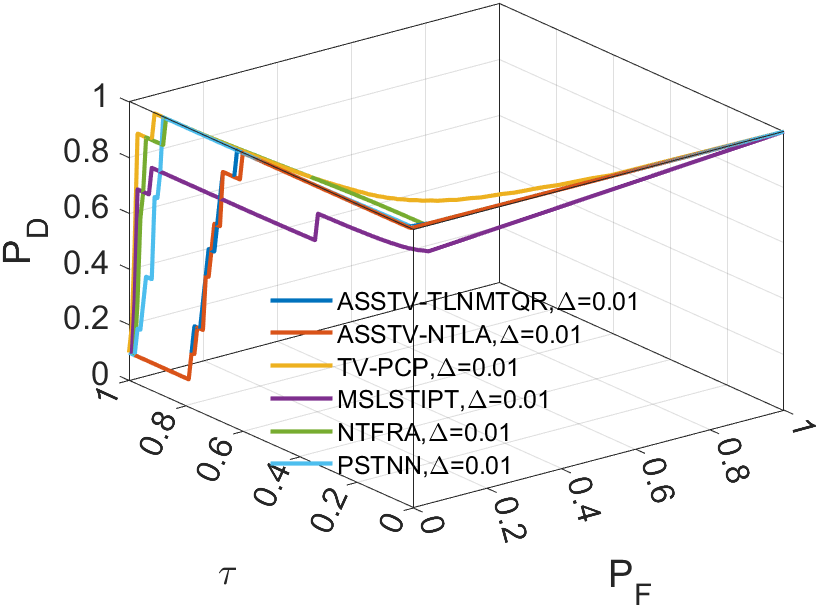}
		\label{sc4-1}
	\end{minipage}
	\begin{minipage}{0.3\textwidth}
		\centering
		\includegraphics[width=0.9\linewidth]{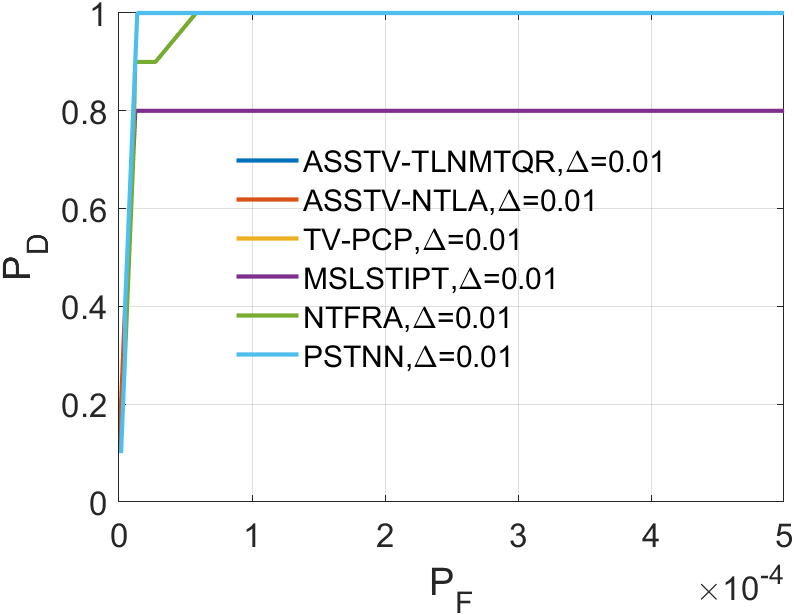}
		\label{sc4-2}
	\end{minipage}
	\begin{minipage}{0.3\textwidth}
		\centering
		\includegraphics[width=0.9\linewidth]{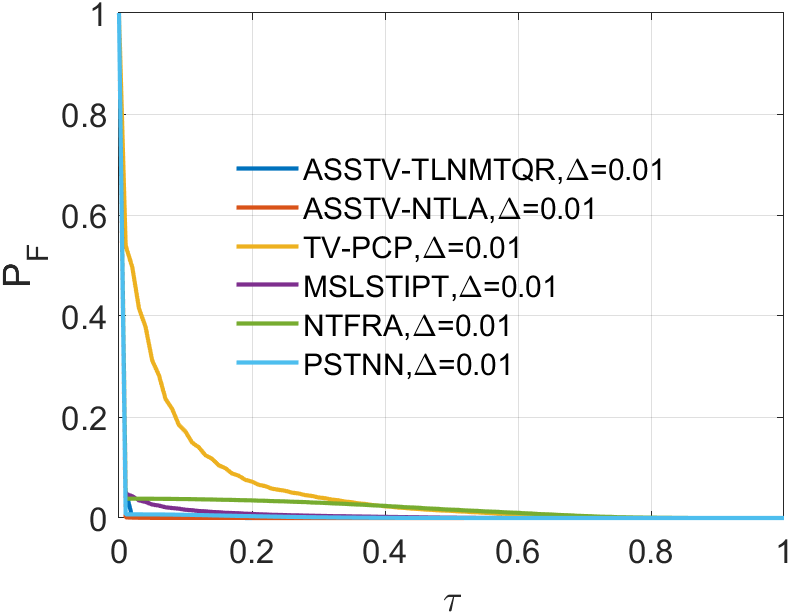}
		\label{4-3}
	\end{minipage}
\end{figure*}
\begin{figure*}[!htb]
	\caption{The 3D ROC of Sequence 5}
	\centering
	\label{sc5}
	\begin{minipage}{0.3\textwidth}
		\centering
		\includegraphics[width=0.9\linewidth]{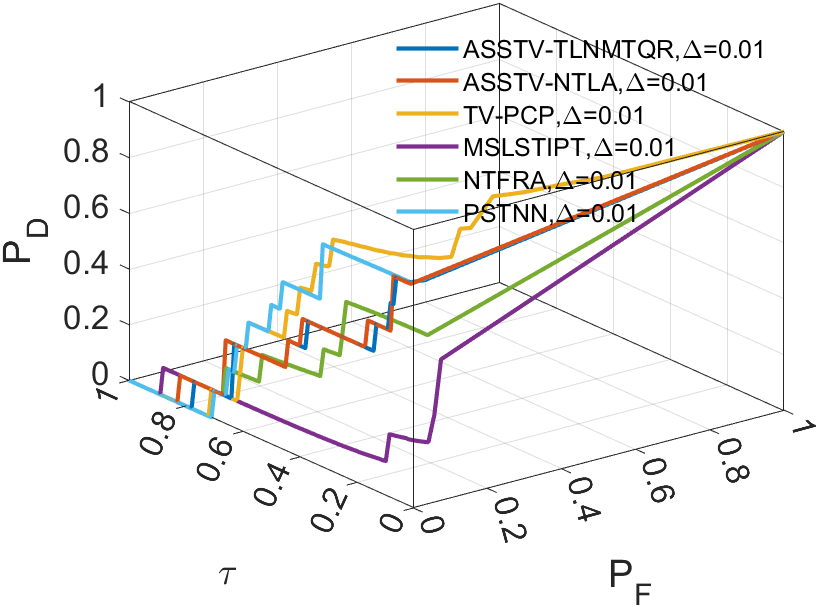}
		\label{5-1}
	\end{minipage}
	\begin{minipage}{0.3\textwidth}
		\centering
		\includegraphics[width=0.9\linewidth]{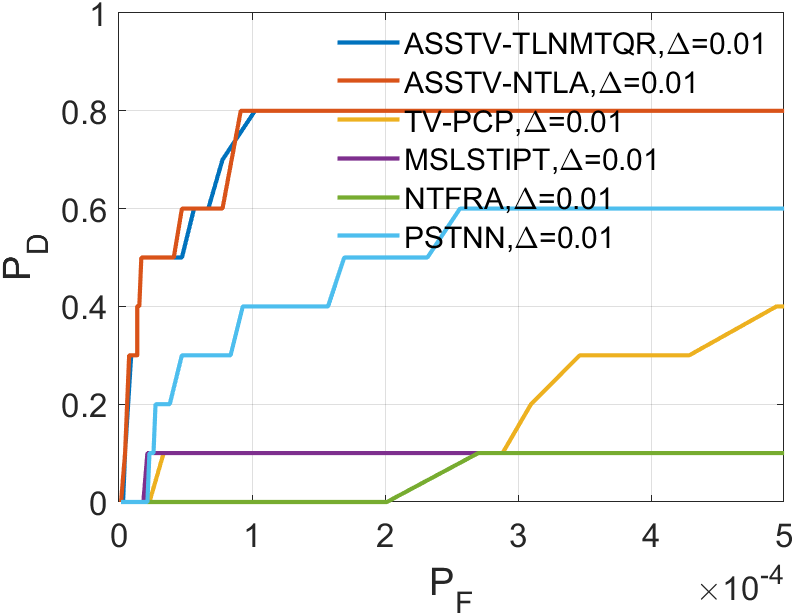}
		\label{sc5-2}
	\end{minipage}
	\begin{minipage}{0.3\textwidth}
		\centering
		\includegraphics[width=0.9\linewidth]{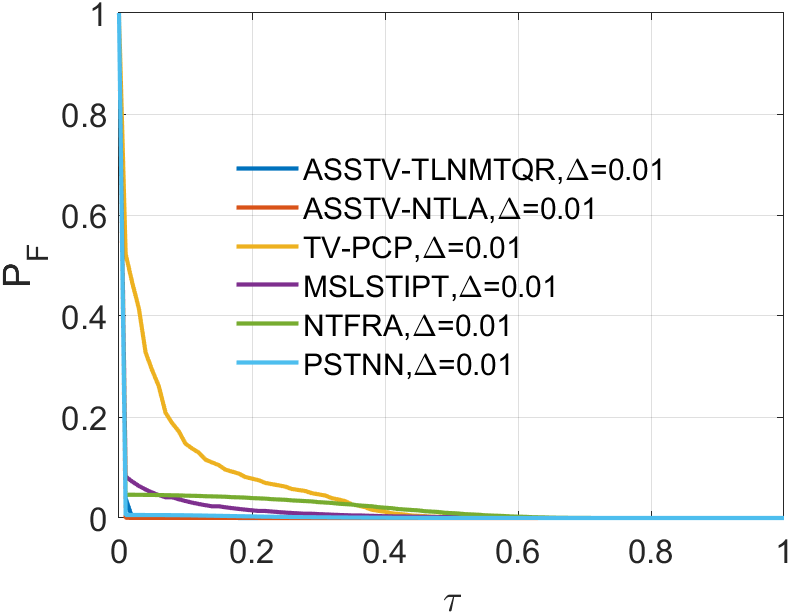}
		\label{sc5-3}
	\end{minipage}
\end{figure*}
\begin{figure*}[!htb]
	\caption{The 3D ROC of Sequence 6}
	\centering
	\label{sc6}
	\begin{minipage}{0.3\textwidth}
		\centering
		\includegraphics[width=0.9\linewidth]{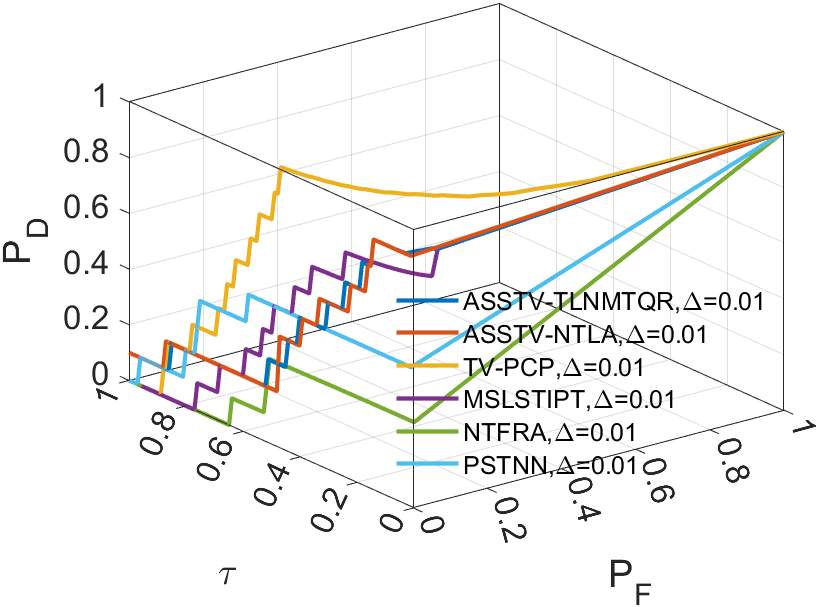}
		\label{sc6-1}
	\end{minipage}
	\begin{minipage}{0.3\textwidth}
		\centering
		\includegraphics[width=0.9\linewidth]{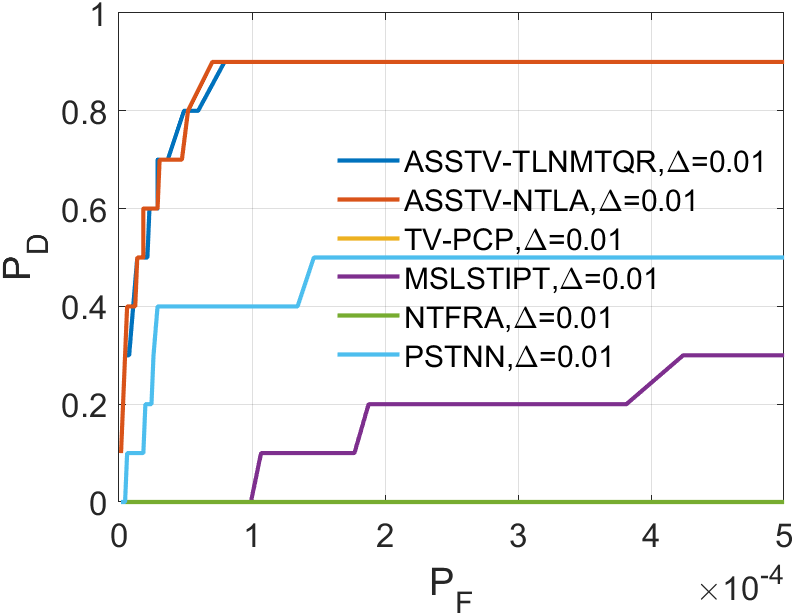}
		\label{sc6-2}
	\end{minipage}
	\begin{minipage}{0.3\textwidth}
		\centering
		\includegraphics[width=0.9\linewidth]{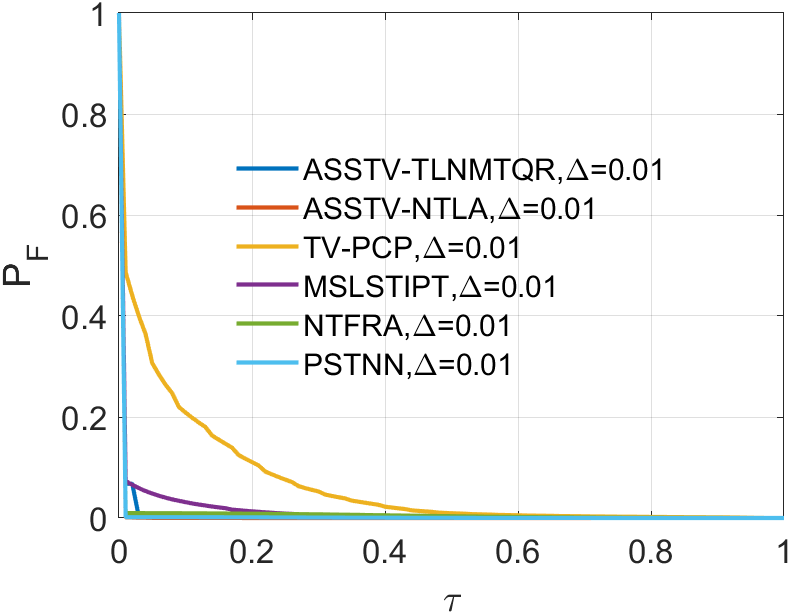}
		\label{sc6-3}
	\end{minipage}
\end{figure*}

\section{CONCLUSION}
This article introduces an innovative tensor recovery algorithm. It combines ASSTV regularization and TLNMTQR methods, resulting in a significant improvement in computational speed. We conducted experiments in infrared small target detection to evaluate the capabilities of our algorithm. In experiments on small infrared target detection, we tested the algorithm's ability to extract the target tensor $\mathcal{T}$ and suppress the background tensor $\cal B$. Our algorithm greatly enhances computational speed without sacrificing accuracy.

In future research, we have several areas for optimization. Firstly, when addressing the $L_{2,1}$ norm minimization problem, we can consider assigning different weights to different eigenvalues to enhance background extraction capabilities. Secondly, apart from parameter adjustments, we can explore more effective methods to determine the optimal size for approximating tensor $\mathcal{K}$.

\small
\bibliographystyle{IEEEtran}\small
\bibliography{IEEEabrv,mybib}

\begin{thebibliography}{10}
\providecommand{\url}[1]{#1}
\csname url@samestyle\endcsname
\providecommand{\newblock}{\relax}
\providecommand{\bibinfo}[2]{#2}
\providecommand{\BIBentrySTDinterwordspacing}{\spaceskip=0pt\relax}
\providecommand{\BIBentryALTinterwordstretchfactor}{4}
\providecommand{\BIBentryALTinterwordspacing}{\spaceskip=\fontdimen2\font plus
\BIBentryALTinterwordstretchfactor\fontdimen3\font minus
  \fontdimen4\font\relax}
\providecommand{\BIBforeignlanguage}[2]{{%
\expandafter\ifx\csname l@#1\endcsname\relax
\typeout{** WARNING: IEEEtran.bst: No hyphenation pattern has been}%
\typeout{** loaded for the language `#1'. Using the pattern for}%
\typeout{** the default language instead.}%
\else
\language=\csname l@#1\endcsname
\fi
#2}}
\providecommand{\BIBdecl}{\relax}
\BIBdecl

\bibitem{signals_matrix}
R.~R. Nadakuditi, ``Optshrink: An algorithm for improved low-rank signal matrix
  denoising by optimal, data-driven singular value shrinkage,'' \emph{IEEE
  Transactions on Information Theory}, vol.~60, no.~5, pp. 3002--3018, 2014.

\bibitem{ISTD1}
Z.~Zhang, C.~Ding, Z.~Gao, and C.~Xie, ``Anlpt: Self-adaptive and non-local
  patch-tensor model for infrared small target detection,'' \emph{Remote
  Sensing}, vol.~15, no.~4, p. 1021, 2023.

\bibitem{ISTD2}
H.~Yi, C.~Yang, R.~Qie, J.~Liao, F.~Wu, T.~Pu, and Z.~Peng, ``Spatial-temporal
  tensor ring norm regularization for infrared small target detection,''
  \emph{IEEE Geoscience and Remote Sensing Letters}, vol.~20, pp. 1--5, 2023.

\bibitem{ISTD3}
L.~Chuntong and W.~Hao, ``Research on infrared dim and small target detection
  algorithm based on low-rank tensor recovery,'' \emph{Journal of Systems
  Engineering and Electronics}, 2023.

\bibitem{compress_sensing}
M.~Lustig, D.~L. Donoho, J.~M. Santos, and J.~M. Pauly, ``Compressed sensing
  mri,'' \emph{IEEE Signal Processing Magazine}, vol.~25, no.~2, pp. 72--82,
  2008.

\bibitem{8260873}
M.~Rani, S.~B. Dhok, and R.~B. Deshmukh, ``A systematic review of compressive
  sensing: Concepts, implementations and applications,'' \emph{IEEE Access},
  vol.~6, pp. 4875--4894, 2018.

\bibitem{CHANDRASEKARAN20091493}
\BIBentryALTinterwordspacing
V.~Chandrasekaran, S.~Sanghavi, P.~A. Parrilo, and A.~S. Willsky, ``Sparse and
  low-rank matrix decompositions,'' \emph{IFAC Proceedings Volumes}, vol.~42,
  no.~10, pp. 1493--1498, 2009, 15th IFAC Symposium on System Identification.
  [Online]. Available:
  \url{https://www.sciencedirect.com/science/article/pii/S1474667016388632}
\BIBentrySTDinterwordspacing

\bibitem{LRSD1}
D.~Bertsimas, R.~Cory-Wright, and N.~A. Johnson, ``Sparse plus low rank matrix
  decomposition: A discrete optimization approach,'' \emph{Journal of Machine
  Learning Research}, vol.~24, no. 267, pp. 1--51, 2023.

\bibitem{LRSD2}
G.~Daniel, G.~Meirav, O.~Noam, B.-K. Tamar, R.~Dvir, O.~Ricardo, and B.-E.
  Noam, ``Fast and accurate t2 mapping using bloch simulations and low-rank
  plus sparse matrix decomposition,'' \emph{Magnetic Resonance Imaging},
  vol.~98, pp. 66--75, 2023.

\bibitem{lin2010augmented}
Z.~Lin, M.~Chen, and Y.~Ma, ``The augmented lagrange multiplier method for
  exact recovery of corrupted low-rank matrices,'' \emph{arXiv preprint
  arXiv:1009.5055}, 2010.

\bibitem{chen2017denoising}
Y.~Chen, Y.~Guo, Y.~Wang, D.~Wang, C.~Peng, and G.~He, ``Denoising of
  hyperspectral images using nonconvex low rank matrix approximation,''
  \emph{IEEE Transactions on Geoscience and Remote Sensing}, vol.~55, no.~9,
  pp. 5366--5380, 2017.

\bibitem{xie2016hyperspectral}
Y.~Xie, Y.~Qu, D.~Tao, W.~Wu, Q.~Yuan, and W.~Zhang, ``Hyperspectral image
  restoration via iteratively regularized weighted schatten $ p $-norm
  minimization,'' \emph{IEEE Transactions on Geoscience and Remote Sensing},
  vol.~54, no.~8, pp. 4642--4659, 2016.

\bibitem{zhang2013hyperspectral}
H.~Zhang, W.~He, L.~Zhang, H.~Shen, and Q.~Yuan, ``Hyperspectral image
  restoration using low-rank matrix recovery,'' \emph{IEEE transactions on
  geoscience and remote sensing}, vol.~52, no.~8, pp. 4729--4743, 2013.

\bibitem{TRPCA}
J.~Xue, Y.~Zhao, S.~Huang, W.~Liao, J.~C.-W. Chan, and S.~G. Kong, ``Multilayer
  sparsity-based tensor decomposition for low-rank tensor completion,''
  \emph{IEEE Transactions on Neural Networks and Learning Systems}, vol.~33,
  no.~11, pp. 6916--6930, 2022.

\bibitem{TLRSD1}
M.~Wang, D.~Hong, Z.~Han, J.~Li, J.~Yao, L.~Gao, B.~Zhang, and J.~Chanussot,
  ``Tensor decompositions for hyperspectral data processing in remote sensing:
  A comprehensive review,'' \emph{IEEE Geoscience and Remote Sensing Magazine},
  2023.

\bibitem{8606166}
C.~Lu, J.~Feng, Y.~Chen, W.~Liu, Z.~Lin, and S.~Yan, ``Tensor robust principal
  component analysis with a new tensor nuclear norm,'' \emph{IEEE Transactions
  on Pattern Analysis and Machine Intelligence}, vol.~42, no.~4, pp. 925--938,
  2020.

\bibitem{7892843}
W.~He, H.~Zhang, and L.~Zhang, ``Total variation regularized reweighted sparse
  nonnegative matrix factorization for hyperspectral unmixing,'' \emph{IEEE
  Transactions on Geoscience and Remote Sensing}, vol.~55, no.~7, pp.
  3909--3921, 2017.

\bibitem{TV}
M.-D. Iordache, J.~M. Bioucas-Dias, and A.~Plaza, ``Total variation spatial
  regularization for sparse hyperspectral unmixing,'' \emph{IEEE Transactions
  on Geoscience and Remote Sensing}, vol.~50, no.~11, pp. 4484--4502, 2012.

\bibitem{TV2}
V.~Vishnevskiy, T.~Gass, G.~Szekely, C.~Tanner, and O.~Goksel, ``Isotropic
  total variation regularization of displacements in parametric image
  registration,'' \emph{IEEE Transactions on Medical Imaging}, vol.~36, no.~2,
  pp. 385--395, 2017.

\bibitem{sstv}
W.~He, H.~Zhang, H.~Shen, and L.~Zhang, ``Hyperspectral image denoising using
  local low-rank matrix recovery and global spatial–spectral total
  variation,'' \emph{IEEE Journal of Selected Topics in Applied Earth
  Observations and Remote Sensing}, vol.~11, no.~3, pp. 713--729, 2018.

\bibitem{ASSTV}
T.~Liu, J.~Yang, B.~Li, C.~Xiao, Y.~Sun, Y.~Wang, and W.~An, ``Nonconvex tensor
  low-rank approximation for infrared small target detection,'' \emph{IEEE
  Transactions on Geoscience and Remote Sensing}, vol.~60, pp. 1--18, 2022.

\bibitem{TSVD}
S.~Weiland and F.~Van~Belzen, ``Singular value decompositions and low rank
  approximations of tensors,'' \emph{IEEE transactions on signal processing},
  vol.~58, no.~3, pp. 1171--1182, 2009.

\bibitem{ZHENG2021108240}
\BIBentryALTinterwordspacing
Y.~Zheng and A.-B. Xu, ``Tensor completion via tensor qr decomposition and
  l2,1-norm minimization,'' \emph{Signal Processing}, vol. 189, p. 108240,
  2021. [Online]. Available:
  \url{https://www.sciencedirect.com/science/article/pii/S0165168421002772}
\BIBentrySTDinterwordspacing

\bibitem{6595533}
C.~Gao, D.~Meng, Y.~Yang, Y.~Wang, X.~Zhou, and A.~G. Hauptmann, ``Infrared
  patch-image model for small target detection in a single image,'' \emph{IEEE
  Transactions on Image Processing}, vol.~22, no.~12, pp. 4996--5009, 2013.

\bibitem{tvpcp}
\BIBentryALTinterwordspacing
X.~Wang, Z.~Peng, D.~Kong, P.~Zhang, and Y.~He, ``Infrared dim target detection
  based on total variation regularization and principal component pursuit,''
  \emph{Image and Vision Computing}, vol.~63, pp. 1--9, 2017. [Online].
  Available:
  \url{https://www.sciencedirect.com/science/article/pii/S0262885617300756}
\BIBentrySTDinterwordspacing

\bibitem{PSTNN}
L.~Zhang and Z.~Peng, ``Infrared small target detection based on partial sum of
  the tensor nuclear norm,'' \emph{Remote Sensing}, vol.~11, no.~4, p. 382,
  2019.

\bibitem{MSLSTIPT}
Y.~Sun, J.~Yang, and W.~An, ``Infrared dim and small target detection via
  multiple subspace learning and spatial-temporal patch-tensor model,''
  \emph{IEEE Transactions on Geoscience and Remote Sensing}, vol.~59, no.~5,
  pp. 3737--3752, 2020.

\bibitem{NTFRA}
X.~Kong, C.~Yang, S.~Cao, C.~Li, and Z.~Peng, ``Infrared small target detection
  via nonconvex tensor fibered rank approximation,'' \emph{IEEE Transactions on
  Geoscience and Remote Sensing}, vol.~60, pp. 1--21, 2021.

\bibitem{ROC}
C.-I. Chang, ``An effective evaluation tool for hyperspectral target detection:
  3d receiver operating characteristic curve analysis,'' \emph{IEEE
  Transactions on Geoscience and Remote Sensing}, vol.~59, no.~6, pp.
  5131--5153, 2021.

\bibitem{ROC2}
S.~Atapattu, C.~Tellambura, and H.~Jiang, ``Analysis of area under the roc
  curve of energy detection,'' \emph{IEEE Transactions on Wireless
  Communications}, vol.~9, no.~3, pp. 1216--1225, 2010.

\bibitem{2D_ROC}
S.~Manti, M.~K. Svendsen, N.~R. Kn{\o}sgaard, P.~M. Lyngby, and K.~S. Thygesen,
  ``Exploring and machine learning structural instabilities in 2d materials,''
  \emph{npj Computational Materials}, vol.~9, no.~1, p.~33, 2023.

\bibitem{istd_data}
B.~Hui, Z.~Song, H.~Fan, P.~Zhong, W.~Hu, X.~Zhang, J.~Lin, H.~Su, W.~Jin,
  Y.~Zhang, and Y.~Bai, ``{A dataset for infrared image dim-small aircraft
  target detection and tracking under ground / air background},'' Oct. 2019.

\end{thebibliography}

\vfill

\end{document}